\documentclass[13pt]{article}
\usepackage{amsmath}
\usepackage{amssymb}
\usepackage{amsthm}
\usepackage{amscd}
\usepackage{amsfonts}
\usepackage{dsfont}
\usepackage{cases}
\usepackage[applemac]{inputenc}
\usepackage{enumerate}
\usepackage[all]{xy}
\usepackage{esint}
\usepackage{graphicx}
\usepackage{url}
\theoremstyle{plain}

\theoremstyle{definition}

\newcommand{\N}{\mathbb N}

\newcommand{\R}{\mathbb R}

\newcommand{\E}{\mathrm E}

\newcommand{\one}{\mathrm\bf 1}

\newcommand{\G}{\Gamma}
\newcommand{\Om}{\Omega}
\newcommand{\om}{\omega}
\newcommand{\Ss}{\cal S}
\newcommand{\al}{\alpha}
\newcommand{\vfi}{\varphi}

\begin{document}
\title{Equilibrium Equations for Human Populations\\  with Immigration}
\date{}

\author{F. Thomas Bruss \\Universit\'e Libre de Bruxelles}
\maketitle

\begin{abstract} \noindent The objective of this article is to create a framework to study  asymptotic equilibria in populations with immigration, and this with a special focus on human populations. We present a new model, based on Resource Dependent Branching Processes, which is now broad enough  to cope with the goal of finding equilibrium criteria under reasonable hypotheses. Our equations are expressed in terms of natality rates, mean productivity and mean consumption of the home-population and the immigrant population as well as policies of the Society to distribute resources among individuals. We also study the impact of {\it integration} of one sub-population into the other one, and in a third model, the additional influence of an ongoing stream of new immigrants.
Proofs of the results are based on classical limit theorems, on Borel-Cantelli type arguments, on the Theorem of envelopment of Bruss and Duerinckx (2015),  on a maximum inequality of Bruss and Robertson (1991), and on an extension of J.M. Steele (2016) of the latter. 
Conditions for the existence of an equilibrium often prove to be severe, and sometimes surprisingly sensitive. This underlines how demanding the real world of immigration can be for politicians trying to make sound decisions. Our main objective is to provide decision help through insights from an adequate theory. Another objective of the present study is to learn  {\it which} of the possible control measures are best for combining feasibility and efficiency to reach an equilibrium, and to recognise the corresponding steps one has to take towards controls. We also make preliminary suggestions to envisage ways to optimal control. 
 As far as the author is aware, all  results are new.


\medskip
\noindent {\bf AMS 2010 Math. subj. classific.}: 60J85, secondary 49J21.

\medskip\noindent{\bf  Short Running title}: Immigration and Equilibrium
\end{abstract}

\section{Motivation}
During the last few years the picture of human migration rates between countries in the world has dramatically changed. Focussing on immigration, the following can be seen. Countries with a long history of immigration such as Australia and the United States, still have a large percentage of  immigrants but the gradient of change is becoming small compared to countries such as  Austria, Denmark, Germany,   Sweden, and others. Until 2013, Germany for example, had an immigration rate of about one per thousand per year (1/1000)/y,
that is below the rates of several other European countries, whereas in 2015 the German exploded to (25/1000)/y. 

Several countries show much goodwill towards immigrants, in particular for refugees in danger in their home countries. Goodwill alone is not sufficient to deal with  problems which arise if  the number of immigrants increases  quickly.
The evident challenges are to offer accommodation to immigrants, and to find employment or create new jobs for them, but also how to give immigrants a real chance to integrate themselves into their new country. Integration seems only possible if the natality rates of immigrants and host population converge sufficiently quickly to each other since otherwise children of immigrants have no surrounding in which they can naturally learn the new home-language. Converging birth rates are thus an important issue. As we shall see,  there are other issues almost as important.

\subsection{Focussing on equilibria}
We will focus our interest on the existence of long-term equilibria. This is a rather theoretical focus, but we think that, provided the models are not unrealistic, the theory we develop is rewarding for a deeper understanding of the effect of immigration.  As a consequence, our results are only indirectly connected with questions surrounding economical and econometric shorter-term aspects of immigration, as e.g. concerning targeted immigration policies, cost-benefit analysis, planning of resource allocation, prediction, and others. For an  authoritative presentation of shorter-term economic aspects of immigration see e.g.  Borjas (2014).

As the present paper will show,
an equilibrium, in fact any kind of equilibrium, may be hard to reach under strong immigration, and it is difficult to give political advice for viable equilibria.  This is our goal, however, and we want to make things as transparent as possible. 

Understanding  the conditions under which a home-population and immigrants can, in the longer run,  attain an equilibrium is a strongly motivated objective. Studies on peace suggest that there can exist no convincing long-term alternative to equilibria.
We shall try in this paper to cope with this challenge by studying a model based on {\it Resource-Dependent Branching Processes} (RDBPs). 

\subsection{Related work} Branching processes with immigration have been studied by several authors, and most of these authors study Markov processes, in particular modified Galton-Watson processes. Classical references are the book by Haccou et al. (2005) and the many articles cited in there. See also the new book on controlled branching processes by Gonz\'ales 
et al.\,(2018). 

As in several other branching process models, a density-dependent development of the population around the so-called critical case will have a natural appeal in our model, and thus our approach shares in part the motivation of the work of Afanasev et al. (2005),  Jagers and Klebaner (2000), Isp\'any (2016), and others. 
Fluctuations under immigration (see e.g. Ispany et. al (2005) and Wei  and Winnicki (1989)) are also naturally at stake although we only speak indirectly about fluctuations.  Concerning density controls,  one would also like to know when, and in what way, a population will leave the region around criticality without control. In this respect our interests come close to those motivating former studies of Bingham and Doney (1974),  Keller et al. (1987), Klebaner and Zeitouni (1994), and, more recently,  Barbour et al. (2015),  Kersting (2018), and Bansaye et. al (2018). 

We cannot  directly profit from these  results because resource dependent branching
processes have a different structure, and our approach must go different ways. It is the notion of society control, which we are forced to incorporate in a realistic view of human behaviour, which explains one part of the major difference in structure. A second one stems from allowing for independence {\it within} each sub-population but sacrificing the independence of sub-processes as such in favour of a {\it common resource space} on which sub-populations have to live. 

When thinking about how to tailor a tractable model, indirect influences can be almost as beneficial as a direct influence. The author sincerely acknowledges  what he has learned from the papers cited above, and  from many others  not mentioned here in the longer branching-process history. All have helped to develop intuition. 
\section {Content and objective of this paper}
 In Section 3 we summarise the notion of RDBPs without immigration, the idea behind them,
 and the reasons why we believe that these processes are an adequate approach  for describing the development of human populations. For the present paper we always understand an RDBP as defined in Bruss and Duerinckx (2015, section 2), and the so-called {\it society's obligation principle} as defined in Bruss (2016, subsection 7.1.1).

In Section 4  we explain why we should not try to use directly the model of Bruss and Duerinckx (2015) if we allow for immigration. Unlike emigration, immigration is indeed not incorporated in their model. The conclusion is that a better model must be found and that we should study a suitable multi-variate process. Although this seems like a natural step, it is less obvious how to do this if we want to study the development of a society which is 
{\it non-discriminating} in the sense that individuals stemming from the home-population or from the immigrant-population are submitted to the {\it same} rules for receiving resources from a {\it common} resource space. Individuals submit random {\it claims} to the society, and for the present paper it suffices to understand a claim as what an individual would like to obtain for individual consumption. The society decides whether to accept or not to accept the claim from the combined (merged) list of claims according to the currently fixed rules conditioned on available resources. This implies that the sub-populations are dependent on each other.

To cope with this problem we recall in Section 5 an extension of what Steele (2016) calls the {\it Bruss and Robertson-inequality} (BR-inequality.) Steele's extension will play a central role since it allows to make full use of the upper bound of the Theorem of envelopment (Th. 4.13 in Bruss and Duerinckx (2015)). This will yield  in the following sections conditions for the survival of both sub-processes. 

 Section 6 defines the notion of an equilibrium. Then it studies a bi-variate RDBP counting individuals from the home-population and  immigrant population without new immigrants as if they behaved like cohabitating resource-dependent populations. This is the adequate model for the case where from some finite time onwards there are no new immigrants and where we want to understand how both  sub-populations would develop. We derive the corresponding criterion for possible equilibria. Here we see that in the so-to-speak typical case, i.e. immigrants are poorer but have more children, an equilibrium cannot be reached without control.
 
 Section 7 incorporates  the feature of {\it integration} where individuals from the immigrant-population become successively part of the home-population and then behave exactly like individuals of the latter. The situation becomes now quite different in the sense that if the society were in complete command of integration, even the unfavourable  typical case may allow for an equilibrium. Increasing the facilities of  integration is indeed one of the easier controls to reach an equilibrium.
 
 In Section 8 we complete the set of our basic models by allowing also an ongoing stream of new immigrants. In particular we can show that our approach stays coherent under a reasonable condition and that the computation of the
 limiting equilibrium follows then the same lines. The important benefit is that the studied different influences can now directly be compared with each other.
  
In Section 9 we examine aspects of the flexibility of our models. We also show why our results, so far all obtained for the so-called {\it weakest-first policy}, are of general interest because many different policies can be re-interpreted as such a policy under  transformed claim distribution functions. Claim distributions can be transformed by mixing different distributions and we return here to the
BR-inequality to understand and interpret its stability with respect to different kinds of mixing. We also shortly address the possibility to modify rules of attribution of resources and to envisage optimal control seen as a problem of optimal transport.
 
 Section 10, finally, collects items of important criticism one may see for our approach and for our models, and then draws the main conclusions.
\section{Resource Dependent Branching processes}

Resource Dependent Branching Processes,  introduced by the author already in 1982 (see Bruss 1984a), are not very known. We will review them  briefly by explaining the motivation behind them. This is done in a summarising style with the intention to facilitate the reading of the present paper, and to draw attention  to  these models. For details we refer to sections 1 and 2 of Bruss and Duerinckx (2015).
\subsection{Features} One part of the idea behind RDBPs is that
a suitable model for human populations must offer several basic features: 
Individuals have to eat, to reproduce, and to work, in order to be able to survive. They consume resources, may inherit and/or save them,  and then again they create new resources for the descendants. 
In our definition of RDBPs we suppose that what is left after consumption will go into a common  {\it resource space,} although many other assumptions would be compatible with the model. And then we need a notion of a society which defines the policy of distributing resources, and also a notion of protest against decisions of the society.

Individual requests (needs) of resources are seen as random variables, called {\it claims}. We can see in our model individual claims as individual consumption, although they are more precisely defined as what individuals request to have at their disposal. If an individual does not receive its claim under the current policy we suppose it shows its protest against Society by emigration before leaving offspring, or equivalently, by not reproducing in the present population.  

Unless stated otherwise a claim is either served completely, or not at all. Such a claim is then consumed or partially consumed, and what is left goes as {\it heritage} into the common reserve for the next generation. Heritage is also seen as resource creation for the next generation. For simplicity, all resulting claims within the same generation are supposed to be 
independent identically distributed (i.i.d.) random variables according to a continuous distribution function $F.$ 
 
Reproduction of human beings is in reality bi-sexual, of course. We refer to Daley (1968), and,  for an overview of such models, to Molina (2010).
 However, since we only study long-term developments for large population sizes it suffices to study reproduction in terms of the {\it average reproduction  of mating units}  (Bruss 1984b, p. 916).  This allows us to always argue as if we had  asexual reproduction, and the mean of reproduction $1$ is understood as the mean reproduction  of mating couples.
 
 The rules to distribute available resources according to the resulting claims are defined by  the current {\it policy} of the population.  For the precise definition see Def. 2.1  in Bruss and Duerinckx (2015).  

\medskip
{\bf RDBPs as local models}

\medskip
The second part of the idea behind RDBPs is that they should serve as {\it local} models. It is not realistic to make long-term hypotheses for the development of human populations. {\it Local} means local in time, that is, defined on a short horizon of one or a few generations. 
An advantage of local models is that they can be tailored with simple assumptions. Assuming  that certain random variables associated with individuals from a human population are i.i.d. within a given generation is easier to defend than assuming that these i.i.d.-hypotheses would hold forever. RDBPs are used in a history-driven set-up to form a global model. This is explained in the Introduction of Bruss and Duerinckx (2015), and made explicit in terms of the {\it society's obligation principle} (Bruss (2016) subsection 1.1.1). 

\subsection{Concatenating RDBPs to a global model} 

The concatenation is as follows.
At time $t-$ we suppose to know the probabilistic prescription of the current RDBP defined on $[t,t+1[$, but not the (precise) prescription of the future RDBPs.
 We suppose that the population has objectives, recalled below, and that the obligation principle forces it to control for these objectives at each time $t-$ with $t\in \N.$ This is done by two actions. First, by examining the current parameters essential for the development (mean rates of birth, production and consumption), and second, by encouraging certain changes of parameters and by changing the policy of distributing resources (see below). 
 
 \bigskip
\noindent The {\it global objective} was defined in Bruss and Duerinckx (2015) by  two natural hypotheses which we maintain throughout the present paper: The large majority of individuals

\begin{quote}
H1: wants to survive and see a future for the descendants, 

\smallskip
H2:  prefers a higher standard of living to a lower one.

 \end{quote}

\noindent Here the hypothesis H1 is supposed to take priority if H2 becomes incompatible with H1, which is often the case.

The follow-up  is supposed to be ruled by the mentioned society obligation principle to observe H1 and H2 which we make now precise. At each time $t-$ the society checks the following question: If all currently observed parameters and the policy to distribute resources were  to stay the same for all future generations - that is, if the current RDBP would run forever- would it then have a positive probability of surviving forever?  If yes, the society keeps this RDBP or, optionally, replaces it by another RDBP fulfilling this requirement. If not, it controls immediately to obtain as quickly as possible a new RDBP for which this answer would be yes. No minimum survival probability is prescribed, provided that it is {\it strictly} positive.

This principle turns the  development of the population into a history-driven sequence of the local RDBPs. We know no details about the future ones, but if the society obligation principle is always respected we know the objective, the possible range of control actions, and thus the possible range of models respecting H1 and H2. This is why we can focus our interest on knowing under which condition a specific RDBP can survive forever.

\smallskip
Note that, unlike the model of a RDBP which has a well-defined probability prescription, this control approach to form a global model for a human society is no probability model. It contrasts therefore other interesting branching process models involving certain forms of competition for resources, as for instance the branching annihilating random walk studied by Perl et al (2015), or, in the context of varying or random environments, the processes studied by Keller et al. (1987), Kersting (2017) and recently Bansaye et al. (2018), Barczy et al. (2018), and Pap (2018). In fact, not being a probability model, the global model seems to contrast any other existing branching process model. Thinking of our objective  and of the complexity of human populations facing an unknown future, the global model may be more adequate, however.  Also,
in each generation decisions of control are based on studying the current RDBP, that is, on a probability model, so that the control decisions themselves are based on a rigorous setting.

\subsection{RDBPs and the {\it wf}-policy}
Let $(\Gamma_n)_{n=1,2, \cdots}$ be an arbitrary RDBP  with mean reproduction of individuals $m$, average productivity $r$ and resource claim distribution function $F(x):=F_X(x)$. It follows from Bruss and Duerinckx (2015) (p. 336, Prop. 4.3) that $(\Gamma_n)$ will get extinct almost surely if it cannot possibly survive under the so-called weakest-first policy ({\it wf}-policy).  This {\it wf}-policy is the policy to distribute the resources  from the available resource space with priority to those individuals with the smallest claims as long as the current resource space allows for it. Recall that those individuals whose claims are not completely satisfied are supposed not to reproduce in the population. Given the random resource claims $X_1, X_2, \cdots, X_n,$ say from
the $n$ descendants in a given generation with available resource space $s>0,$ the total number of those who will reproduce with the population is thus \begin{align} N(n,s):=\sup \{1\le k \le n:\sum_{j=1}^k X_{j,n}\le s \}, \end{align}
where $X_{j,n}$ is the $j$th smallest order statistic of the $X_1, X_2, \cdots, X_n.$ The random variable $N(n,s)$ is thus the counting variable for the {\it wf}-policy. It is easier to understand and to deal with than counting variables of arbitrary policies, and this is one reason why it plays a major rule throughout this paper. 

The second reason is that, as shown in Section 9, the {\it wf}-policy can be adapted to quite a large class of policies.
Here already the essence of the idea:
Suppose  that the society decides to quit the {\it wf}-policy by, for instance, ignoring all claims falling in certain subintervals of the positive half-line, or else, accepting claims in some other subintervals with some higher probability. If the society announces this change of policy at the beginning of a generation, then individuals may reconsider their claims, and the original claim distribution function $F$ is likely to change in the next generation into some
other distribution function $\tilde F.$ Although the shift of claims may be  difficult to predict in practice, a {\it wf}-policy with respect to $\tilde F$ is in general different from the {\it wf}-policy with respect to $F,$ that is,  the society
applies now, in terms of $F$, another policy.

The class of policies which, for a given $F,$ can be presented as a {\it wf}-policy under some modified distribution $\tilde F$ can be shown to be comfortably large because  a subclass of this class, tentatively called  {\it pure-order policies} by Bruss and Duerinckx (work in progress) is  large enough for most practical purposes. Actually, for the essence of our objective in this paper the mentioned idea of relocating claims will be sufficient. 
This is why we will confine our interest in the present paper until Section 8 included to the {\it wf}-policy. 

\section{Subtleties in understanding immigration}
The subtlety in understanding immigration is best visualised by looking first at the original (univariate) RDBP  introduced in Bruss and Duerinckx (2015). We recall that $m$ denotes the offspring mean of an individual, and $r$ its average production of resources. 
\subsection{Scarce resources}
Confining to the economically relevant case of {\it scarce resources,} we suppose that the average amount of resources $r$ left by an ancestor does not exceed the average total sum of claims submitted by his descendants.  With $m$ denoting the reproduction mean  of an individual, and $F$ the distribution function of individual claim sizes with mean $\mu$, this condition translates into
\begin{align}r \le m \int_0^\infty xdF(x)=m\mu.\end{align}
As recalled before, the {\it wf}-process can only survive if $mF(\tau)\ge1,$ where the parameter $\tau:=\tau (F,r,m)$ is defined by\begin{align}\tau=\inf \left\{t\ge 0: \int_0^t x dF(x) \ge r/m.\right\}\end{align} If $F$ is seen as being a fixed claim distribution, we can drop $F$ and write \begin{align}\tau:=\tau(r,m):=\tau(F,r,m).\end {align} We think of claims
as being evaluated in monetary units, and we assume  that $F$ is strictly increasing and absolute continuous in some neighbourhood of this solution. Hence $\tau$ is uniquely determined by \begin{align} \E\left( X \,\one_{\{X\le \tau\}}\right) = \int_0^\tau x dF(x) = r/m.\end{align} 

 The product $mF(\tau)$ can be seen as the {\it effective long-run multiplication rate} of the process when all other factors are kept invariant. The  condition $mF(\tau)\ge1$ reminds us of the criticality/super-criticality condition for a Galton-Watson branching process. RDBPs are much more complicated processes, of course, but the comparison is partially justified in as much as those individuals which reproduce within a given generation  do so independently of each other.
 
The integral equation (3) yielding the long-run multiplication rate $mF(\tau)$ yields  more by studying the gradient of change  of $mF(\tau(m,r))$ if the parameters $r$ and $m$ change and interact. 
 Let us first look at the influence of the parameter $r$ (average resource production of an individual) for 
 fixed natality $m.$ It was shown in Bruss (2016) that \begin{align}\frac{\partial \left(m F(\tau(r,m))\right)}{\partial r} \ge 0,\end{align} and that this is strongly related with the extinction probability.
Typically, if $r$ goes down the chance of survival decreases, although this may seem a priori
unrelated. (This is, by the way, a delicate observation for those countries in which the birth rate $m$ alone is already below 1 and which thus will get extinct anyway. Decreasing $r$ accelerates extinction, and since decreasing the age of retirement of individuals reduces their life-productivity, and thus reduces $r,$ early retirement is harmful for the probability of survival.)
The influence of a change of $m$ for fixed $r$ and $F$ hides no surprise. It fits our intuition, namely, if $m$ increases, the survival probability increases  (Bruss (2016 ), Theorem 7.7 (ii)).

\medskip
Now comes the important question 
what will happen if $r$ and $m$ change at the same time?  This is what often occurs with immigration in the real world, because, in general, the poorer populations have higher birth rates and a smaller expected productivity. Hence $m$ goes up and $r$ goes down. As we have just seen above the survival probability can now increase only if the influence of the increase of $m$ on the crucial product $mF(\tau)$ is stronger than the negative influence caused by the decrease of $r.$ 
Since $\tau$ (see (3)) is an implicit function of $m$ and $r$ it is  hard to see what will happen. Moreover we can no longer speak of a long-term multiplication rate $mF(\tau)$ because, a priori,
it is not meaningful to assume that $F$ is fixed. Hence the simplification $\tau(r,m):=\tau(F,r,m)$ is no longer justified. The point is that when a larger number of people with different cultural and economic background joins the home-population, this will change the distribution function $F$ of claims. 

The difficulty induced by immigration is that the mechanism of this change is not at all transparent. It might seem reasonable to push our analysis through by imposing that $F$ belongs to a set of distribution functions in some class parametrised by $m$ and $r,$
$F\in \{F^{(m,r)}; m, r \in \R ^+\}$ , say. However, it is not realistic to suppose that we understand the interaction of the three assumed actors $r, m$, and a probably delayed result $F^{(m,r)},$ and we must attack the problem in a different way.

\smallskip
It is a more recent result of J. M. Steele (2016) extending an inequality of Bruss and Robertson (1991) which instigated the idea of how to do this in a tractable way.
\section{Maximum inequality and Steele's extension} Let $X_1, X_2, \cdots$ be a sequence of positive random variables with respective absolute continuous distribution functions $F_k,~k=1, 2, \cdots,$ and let $n>0$ be a fixed positive integer. Further let $ X_{1,n} \le X_{2,n} \le \cdots \le  X_{n,n}$ be the increasing order statistics of $X_1, X_2, \cdots ,X_n$, and let for $s \in \R^+, $ 
\begin{align}\tilde N(n,s)=\begin{cases} 0, {~\rm if}~ X_{1,n}>s,\\\sup\{k \in \N: X_{1,n}+X_{2,n}+ \cdots +X_{k,n} \le s\}, \,{\rm otherwise.}\end{cases}\end{align} $\tilde N(n,s)$ is thus essentially the same as $ N(n,s)$ defined before  in (1), namely the maximum number of variables in $\{X_1,X_2, \cdots, X_n\}$ we can sum up without exceeding $s$, the only difference being that the $n$ order statistics $X_{1,n}\le X_{2,n}\le \cdots \le X_{n,n}$
are not necessarily the order statistics of identically distributed random variables. In the following we reserve the notation $N(s,n)$ for the case of identically distributed random variables. 

We note that both $N(n,s)$  and $\tilde N(n,s)$ are quasi-stopping times in the (more precise) sense that $N(n,s)+1$  and $\tilde N(n,s)+1$ are stopping times
on the sequence of the corresponding increasing order statistics.
Although we will not directly use this fact, it is helpful for the intuition for the following results.

\bigskip \noindent {\bf Theorem 5.1} {(Bruss and Robertson (1991)) }
\begin{align} 
\E(N(n,s))\le nF(\tau),
\end{align} 
where $\tau:=\tau(n,s)$ solves \begin{align}~n \int_0^\tau x dF(x)= s.\end{align}
If, moreover, the $X_k$'s are {\it independent} and
$(s_n)\to \infty$ with $\lim s_n/n>0$ then \begin{align} n^{-1}N(n,s_n)/F(\tau(n,s_n))~ \to~1~\rm{a.s.}, ~ as ~ n \to \infty.\end{align}

\smallskip \noindent For the proof of the inequality (8) with (9) see Lemma 4.1 of Bruss and Robertson (1991), page 622;  for the proof of (10) see Theorems 2.1 and 2.2 of the same paper.  

\medskip\noindent It is the inequality (8) with (9) which attracts here  our main interest.  Steele (2016) called this inequality the {\it Bruss-Robertson inequality} (BR-inequality). His extension of (8) and (9) displays a more versatile inequality.

\medskip \noindent
{\bf Theorem 5.2} {(J. M. Steele (2016)) } Let $X_1, X_2, \cdots, X_n$ be such that each $X_k$ has a absolute continuous distribution function $F_k,$
and let $\tilde N(n,s)$ be defined as in (7). Then
\begin{align} \E(\tilde N(n,s)) \le \sum_{k=1}^n F_k(\tau), \end{align} where  $\tau:=\tau(n,s)$ is a solution of\begin{align} \sum_{k=1}^n \int_0^\tau x dF_k(x)= s. 
\end{align} For the proof see section 3 in Steele (2016).

\medskip\noindent 
 Hence Steele (2016) drops the assumption that the $X_k$ are identically distributed; each $X_k$ can now have its own continuous distribution $F_k$, and the corresponding result remains true. 
\
Bruss and Robertson (1991) were motivated by problems in which the result (10) played the main role, and where it was natural to suppose the $X_k$'s to be i.i.d. random variables.  Although seeing that independence was not used in their proof of (8) and (9) in  Theorem 5.1 they did not point this out.  Interestingly, Steele's proof is a skilfully adapted version of the proof of Bruss and Robertson. Moreover, as we will see in  Section 9, in some cases there is some benefit in re-interpreting Steele's extension as the BR-inequality, or vice-versa, for  mixed distributions.

 \medskip 
It is Steele's great merit (Steele (2016)) to have underlined  the true interest of this result without the assumption of identically distributed $X_k$'s, and that it is not the joint distribution of variables which counts but only their marginals. Steele's extension was an eye-opener for the author and instigated the author's approach presented in this paper. It intervenes repeatedly in the important proofs. 

We should also mention here that Steele (2016) gives examples of applications strongly related with the work of Samuels and Steele (1981), Arlotto et. al (2015),  and Bruss and Delbaen (2001) in the  domain of monotone subsequence problems, but also examples hinting to quite different problems, as e.g. in combinatorial problems.  Again differently motivated, they are of independent interest, and many readers may find them very stimulating.

\section {New RDBP-model}

 We are now ready to study cohabitation of sub-populations, and also immigration. 
 
 \smallskip The idea is to use Steele's extension (Theorem 5.2) for different classes of parameters and  different claim distribution functions. 
By different classes we mean essentially two, namely those associated with the home-population, $C_h$, say, respectively the immigrant-population, $C_i$, say. In Section 8 we will also refer to an additional class of new immigrants $C_{ni}$. In all notations used in the present paper the indices $h$, $i$, and $ni$ are mnemonic for  home-population, immigrant-population, and new immigrants, respectively.  

\smallskip
We first define the basic RDBP with immigration in terms of two sub-populations living under a common constraint of resources. 
\medskip

\noindent{\bf Definition 6.1}:
Let $(\Gamma(t))_{t=0,1,2, \cdots}:=(\Gamma^h_t,\Gamma^i_t)_{t=0,1,2, \cdots}$ be a bivariate counting process with values in $\N^2$ defined on a filtered probability space $(\Om,{ \Ss}, ({\Ss}_t)_{t \in \N}, P),$ where
\begin{quote}
\smallskip
(i) $\Gamma^h_0=h_0, \Gamma^i_0=i_0,$  with $h_0, i_0 \in \{1, 2, \cdots\},$

\smallskip
(ii) ${(\Gamma^h_t)_{t=0, 1, \cdots}}$ is a RDBP with mean reproduction (mean number of descendants) $m_h$, a mean resource space contribution $r_h,$ an individual claim size distribution $F_h.$ The corresponding mean claim is denoted by $\mu_h=\int_0^\infty x dF_h(x),$

\smallskip
(iii) ${(\Gamma^i_t)_{t=0, 1, \cdots}}$ is a RDBP with corresponding parameters $m_i, r_i$  corresponding claim size distribution function $F_i,$ and mean claim $\mu_i=\int_0^\infty x dF_i(x),$

\smallskip
(iv) $(\Gamma^h_t)$ and $(\Gamma^i_t)$ are supposed to be submitted to the same policy of resource distribution from the common resource space built up by the sum of all individual resource contributions provided by the bi-variate process $(\Gamma(t))_{t=0,1,2, \cdots}.$ \end{quote}

Here it is understood that, whenever we speak of a RDBP, all assumptions of Bruss and Duerinckx (2015) are supposed to be satisfied. We also recall that, in order to assure almost-sure convergence of sample means in the rows of the arrays of the involved random variables, we sometimes needed complete convergence (see e.g. Asmussen and Kurtz (1980)), and this is why we suppose that all second moments of the random variables exist. 
As  in most branching processes of interest, we suppose that the probability of an individual having no offspring is strictly positive for both sub-populations.
We maintain these hypotheses throughout this paper.

Note that, in using RDBPs to model the sub-processes $(\Gamma^h_t)$ and $(\Gamma^i_t)$, reproduction and resource space contributions of individuals are i.i.d, random variables within each sub-process separately. This is intrinsic in the definition of an RDBP. The processes $(\Gamma^h_t)$ and $(\Gamma^i_t)$ are however, without further assumptions, not independent of each other, because of (iv).  In analogy to the case of {\it scarce resources} for one population, (see (2)), we will assume, in all what follows, that the expected total production of resources of all sub-populations together is less than the expected sum of all their claims together.

Also, since we have so far no ongoing flow of new immigrants joining the home-population, and no integration of one sub-population into the other one, we see the model defined by (i)-(iv) as a model of {\it cohabitation}. We refer to it as Model I.

 \subsection{Equilibria for Model I}
 
 \smallskip
Consider in Model I a fixed generation $t>0$ and the transition from state $(\G^h_t,\G^i_t)$ to state $(\G^h_{t+1},\G^i_{t+1}), t=1, 2, \cdots.$ Let $D_t^h(u)$ (respectively, $R_t^h(u)$) be the random number of offspring (respectively, random total resource contribution) of $u$ individuals of the home-population in generation t, and let $D_t^i(u)$ and $R_t^i(u)$ be defined correspondingly for the immigrant population.  Given $\G_t^h$ and $\G_t^i$ the total resource space created by the two together equals, according to (iv),  $\tilde R(\G_t^h,\G_t^i):=R^h_t(\G_t^h)+R^i_t(\G^i_t).$

 For the {\it wf}-policy applied to the joined population we have from Steele's extension (see (12)) the corresponding random equation
\begin{align} 
D_t^h(\G^h_t)\int_0^{\tau_t}xdF_h(x)+
D_t^i(\G^i_t)\int_0^{\tau_t}xdF_i(x)=\tilde R(\G_t^h,\G_t^i),\end{align}
where $\tau_t$ is also a random variable, namely according to $\tau(n,s)$ defined in Theorem 5.2,  $$\tau_t:=\tau_t\left(D_t^h(\G^h_t)+D_t^i(\G^i_t), R^h_t(\G_t^h)+R^i_t(\G^i_t)\right).$$
Note that both random equations are well-defined for all $t=1, 2, \cdots$ and all $\om\in \Om$ with the distribution functions of random claims $F_h$ and $F_i,$ as before, not depending on $\om \in \Om.$ 
We now define first the notion of an equilibrium.

\bigskip
\noindent{\bf Definition 6.2} We say that the bivariate process $(\G(t))_{t=1, 2, \cdots}$ tends to an {\it equilibrium}, if there exists a random variable $\alpha$ defined on $(\Om,{ \Ss}, ({\Ss}_t)_{t \in \N}, P),$ with $0< \alpha < \infty$  such that  \begin{align}P\left( \lim_{t\to \infty} \frac{\G^i_t}{\G^h_t} =\alpha \,\Big| \,\G^i_t\not \to 0, \G^h_t \not \to 0\right)=1.\end{align}
{\bf Remark 6.2}: We thus understand an equilibrium as a non-trivial equilibrium between the two sub-populations, that is we do not include $\alpha=0$ or $\alpha^{-1}=0$ in the definition. We will see in Subsection 6.5.2  that the set of possible values of $\alpha$ (seen as "candidates values" for an equilibrium) is typically very small, and in realistic situations often consisting of at most one point. 

\subsection{Conditions for an asymptotic equilibrium without new immigrants}

Recall the Envelopment Theorem (see p. 314, Th. 4.14, Bruss and Duerinckx 2015). Its last statement says that if a given RDBP dies out with probability one under the {\it wf}-policy (written as $q_W=1$) then any other RDBP with the same parameters and claim distribution  would die out with probability one ($ q_\G=1$ for all $\G$).
 RDBP's with the same parameters and claim distribution can only differ in their policies. Hence, in other words, no change of policy whatsoever can enable a process to survive with a positive probability, if the corresponding {\it wf}-process dies out with probability one. 
 Since survival is necessary for the existence of a 
(nontrivial) equilibrium, this is a central result in what follows. Moreover, as said before, studying our process $(\G_t)$ under this specific {\it wf}-policy is less restrictive than it may look.

\medskip Before stating the first main result, a remark on notation.  The existence of the random variable $\alpha$ in Definition 6.2 necessitates of course the existence of candidate values $\alpha_1, \alpha_2,\cdots.$ For easy of notation we  use, whenever this leads to no ambiguity, the notation $\alpha$ for a (fixed) candidate value.
 
 \bigskip
\noindent {\bf Theorem 6.1}

\smallskip\noindent
(a) Let $\cal S$ be the union of the supports of the claim size distributions $F_h$ and $F_i.$ If the natality means $m_h$ and $m_i$ as well as the productivity means $r_h$ and $r_i$ stay invariant over all generations, then an equilibrium can only exist in Model I if there exists a value $0<\alpha<\infty$ and a corresponding value $\tau:=\tau(\alpha) \in \cal S$  satisfying the equation
\begin{align}m_h\int_0^\tau xdF_h(x) + \alpha\,m_i\int_0^\tau xdF_i(x) =r_h +\alpha r_i.\end{align}
subject to the constraints
\begin{align}m_h\, F_h(\tau)=m_i\, F_i(\tau) \ge1.\end{align} (b) Moreover, conditioned on the event that $\{\G_t^i/\G_t^h\not\to 0\}\cap \{\G_t^i/\G_t^h \not\to \infty\}$ replacing the constraints
(16) by $m_h\, F_h(\tau)=m_i\, F_i(\tau) >1$ implies that (a) becomes also a sufficient condition for the existence of an equilibrium.

\bigskip
 \noindent{\bf Proof}: 
The proof consists of four parts, the first three  (i)-(iii) proving (a), and (iv) proving (b).

\medskip
\noindent (i)~We will  first show that if a limiting equilibrium exists then necessarily both sub-processes tend to infinity as $t\to \infty$, that is,\begin{align} \lim _{t \to \infty}\left(\G^i_t/\G^h_t \right){\rm exists~}\implies  P(\G^i_t \to \infty)|\G^i_t \not \to 0)\,=\,P(\G^h_t \to \infty| \G^h_t \not \to 0)=1.\end{align}

\smallskip
The proof of part (i) is by contradiction. 

Suppose that the statement (17) is wrong. We first note that both processes $(\G^i_t)$ and $(\G^i_t)$ must stay bounded away from zero as $t\to\infty,$ since, by definition, zero is an absorbing state for both processes because there are no new immigrants after time $0.$ Since $\alpha$ must satisfy $0<\alpha<\infty$, no sub-population may disappear. But then, if (17) is false, this means that there exist bounds $b_h>0$ and $b_i>0$, say, such that $$P(\G^i_t \le b_i ~\rm{i.o.}| \G^i_t > 0) >0 {~\rm and/or~}P(\G^h_t \le b_h ~\rm{i.o.}| \G^h_t > 0) >0,$$ where i.o. stands for {\it infinitely often.} Put $b=\max\{b_i, b_h\}$ and $p_0=\min\{p_0^i, p_0^h\},$  where $p^h_0>0$, respectively $p^i_0>0$, denotes the probability, that a randomly chosen individual in the home-population, respectively immigrant-population, has no offspring. Since reproduction of individuals is mutually independent {\it within} each sub-population, we must have $$\sum_t
P(\G^h_{t+1}=0|\G^h_t>0) =\infty {~~\rm or~~}\sum_t
P(\G^i_{t+1}=0|\G^i_t>0) = \infty$$ because in both sums all terms are non-negative, and, in at least one sum, infinitely many terms are greater than or equal $p_0^b>0. $ This implies  (see e.g. Corollary 1 of Bruss (1980)) that at least one sub-process will get extinct almost surely. This is in contradiction to (14), however, and hence, conditioned on survival of both sub-processes,
$$\G^i_t \to \infty {\rm~a.s.} \hbox{~and}~ \G^h_t \to \infty{\rm~a.s.~}{,\rm as~} t\to \infty,$$ as stated in (17).

\bigskip
\noindent (ii)~We now prove that, for a given $0<\alpha<\infty$ satisfying (14), there must exist a value $\tau=\tau(\alpha)$ such that equation (15) is satisfied.

First note that if such a value $\tau$ exists for a given value $\alpha,$ then $\tau$ is unique if  the densities $dF_h(t)/dt$ and $dF_i(t)/dt$ do not vanish at the same time in a neighbourhood of $\tau$ because both integrands on the l.h.s. of equation (15) are non-negative. If we denote $\cal S$ the union of the supports of $F_h$ and $F_i$ we can define more generally $$\tau:=\inf\left\{t\in {\cal S}: m_h\int_0^t xdF_h(x) + \alpha\,m_i\int_0^t xdF_i(x) =r_h +\alpha r_i\right\}.$$ This implies the uniqueness of $\tau$ in any case, and justifies the notation $\tau:=\tau(\alpha).$

\smallskip
We now turn to equation (13) with  $\tilde R(\G_t^h,  \G^i_t) = R^h_t(\G_t^h)+R^i_t(\G^i_t).$

Since reproduction and resource production of individuals are independent variables within each sub-process, and since  $F_h$ and $F_i$ are fixed distribution functions, we can apply the strong law of large numbers in equation  (13) for both processes separately.  Moreover, we will see at the same time that $\tau_t$ converges to a constant $\tau:=\tau(\alpha)$ almost surely. 

Indeed, by dividing both sides of (13) by $\G^h_t,$  and using the dummy multiplication factor $\G_t^i/\G_t^i$ for the second terms on both sides, this equation becomes\smallskip
\begin{align} 
\frac{D_t^h(\G^h_t)}{\G^h_t}\int_0^{\tau_t}xdF_h(x)+
\frac{D_t^i(\G^i_t)}{\G^i_t}\frac{\G^i_t}{\G^h_t}\int_0^{\tau_t}xdF_i(x)=\frac{R_t^h(\G_t^h)}{ \G_t^h}+\frac{R^i_t(\G^i_t)}{\G^i_t}\frac{\G^i_t}{\G_t^h}.\end{align}
Accordingly, conditioned on survival of both sub-processes $(\G_t^h)$ and $(\G_t^i)$, the term multiplying the first integral on the l.h.s. converges almost surely to $m_h$  whereas the first term on the r.h.s. almost surely to $r_h.$  Moreover, if the limit $\alpha$ in (14) exists then, as $t \to \infty,$
$$\frac{D_t^i(\G^i_t)}{\G^i_t}\frac{\G^i_t}{\G^h_t}\to \alpha\, m_i {\rm~~ a.s.}~{\rm ~and}~\frac{R^i_t(\G^i_t)}{\G^i_t}\frac{\G^i_t}{\G_t^h}\to \alpha\, r_i {\rm~~ a.s.} $$ Hence, conditioned on survival of both sub-populations and on the existence of the limit $\alpha$, the r.h.s. of (18) has the limit $r_h+\alpha r_i$ a.s. as $t \to \infty.$ This implies that the corresponding l.h.s. must  also have  a  limit.  Since the upper bound is the same $\tau_t$ in both integrals of equation (18), and both integrals have non-negative integrands, we conclude  that $\tau_t$ must converge almost surely  to a constant $\tau$, as $t\to \infty$ . 
Taking these two arguments together we conclude that, if an equilibrium exists, then the  corresponding $\alpha$ and $\tau$ must satisfy the limiting analogue of (18), namely 
$$m_h\int_0^\tau xdF_h(x) + \alpha\,m_i\int_0^\tau xdF_i(x) =r_h +\alpha r_i.$$ This is equation (16) as claimed in the Theorem.
Moreover, with our definition of $\tau$ for a given $\alpha$, the $\tau:=\tau(\alpha)$  must  be unique. 
This proves part (ii).

\bigskip\noindent
(iii) In order to see why the combined constraint qualifications (16) must hold for $\tau=\tau(\alpha),$ we first prove the  equality part of it (which, at first, may look surprising). Let the random variable $\alpha_t$  be defined by $\alpha_t=\G^i_t/\G^h_t,$ or equivalently\begin{align} \frac{1}{1+\alpha_t}=\frac{\G_t^h}{\G_t^h+\G_t^i}, ~ t=1,2 \cdots.\end{align}
Conditioned on survival of both sub-processes we  know thus from part (i) that $\alpha_t\to\alpha$ for some $\alpha$ and, as seen in part (ii), $\tau_t \to \tau:=\tau(\alpha)$ a.s. as $t\to \infty.$ 

Further, in order to be an element of the home-population at time $t+1,$ it is necessary for an individual to be a descendant of it,   and sufficient if its resource claim does not exceed the threshold $\tau_{t}.$ 
It follows that, conditioned on survival of  $(\G_t^h)$ and $(\G_t^i),$ the random fraction of individuals belonging to
the home-population one generation later can therefore be written as
\begin{align}\frac{1}{1+\alpha_{t+1}}=\frac{D_t^h(\G_t^h) F_h(\tau_t)}{D_t^h(\G_t^h) F_h(\tau_t)+D_t^i(\G_t^i) F_i(\tau_t)},  \end{align}
where $\tau_t \to \tau$ a.s. as $t \to \infty.$

Now divide on the r.h.s. of this equation the numerator and denominator by $\G_t^h.$ Using part (i) of the proof and the existence of the limit $\alpha$ we see  then that, conditioned on survival of both sub-proceeses,  $$D^i_t(\G_t^i)/D^h_t(\G_t^h)\to m_i \alpha /m_h  \rm{~a. s. ~as~} t \to \infty.$$ Taking the limit on both sides of (20) for $t\to \infty$ yields then after straightforward computations
\begin{align} \frac{1}{1+\alpha}=\frac{1}{1+\alpha \left(m_i F_i(\tau)/m_hF_h(\tau)\right)},\end{align}
and hence $m_i F_i(\tau)=m_hF_h(\tau).$ 
This proves the equality part of the constraint qualification (15).

\medskip To complete the proof of part (iii) it remains to show that  the conditions $m_hF_h(\tau) \ge 1$ and $m_hF_h(\tau) \ge 1$  are necessary for the existence of an equilibrium.

\smallskip
We argue again by contradiction. 

Suppose the contrary, and suppose first that $0\le m_hF_h(\tau) < 1.$ Let $b \in \,] m_hF_h(\tau), 1[.$  Since  we can confine our interest on the case $\G_t^h \to \infty$ almost surely as $t \to \infty$, and since $D_t^h(\G_t^h)$ is the sum of $\G_t^h$ i.i.d. random variables, we see straightforwardly from the strong law of large numbers and $b<1$ that, for all $\G_t^h$ sufficiently large, $$\E\left(\G_{t+1}^h\big|\G_t^h>0\right) = \E\left(D_t^h(\G_t^h)F_h(\tau_t)\big|\G_t^h>0\right)< \E\left(b\G_t^h\big| \G_t^h>0\right)<\E\left(\G_t^h\big| \G_t^h>0\right).$$ This implies that the process $(\G_t^h)$ stays bounded in conditional expectation (conditioned on non-extinction) so that $$\sum_t (p_o^h) ^{\E(\G_{t+1}^h|\G_t>0)} = \infty,$$ where we recall that $p_o^h$ denotes the probability that an individual in the home-population has no children. But then, by another Borel-Cantelli type argument related with the one we gave before (see now e.g. Bruss (1978), pp 54-56, Theorem 1) we get the contradiction $(\G_t^h) \to 0 {\rm ~a.s.}$ as $t \to \infty.$
A contradiction is obtained in an analogous way by supposing that  an equilibrium exists and $m_iF_i(\tau) < 1.$ 

\smallskip This completes the proof of part (iii).

\bigskip \noindent
(iv) We finally have to show that, given that $\{\G_t^i/\G_t^h\not\to 0\}\cap \{\G_t^i/\G_t^h \not\to \infty\},$ then, replacing the condition
(16) by the slightly stronger condition\begin{align*}m_h\, F_h(\tau)=m_i\, F_i(\tau) >1\end{align*} is  sufficient  for an $\alpha$-equilibrium to exist with a strictly positive probability.  

As we know that the equality
$m_h\, F_h(\tau)=m_i\, F_i(\tau)$  is necessary for the existence of an equilibrium, as we have shown already in the first part of the proof of part (iii), we only have to show to shown that $$P\Big(\{\G_t^h \to \infty\}\cap\{\G_t^i\to \infty\}\Big|m_hF_h(\tau)=m_iF_i(\tau)>1\Big)>0,$$
because, given the joint event $\{\G_t^h \to \infty\}\cap\{\G_t^i\to \infty\},$  we can follow the arguments from (18) up to (21) to establish the limiting equation (15).

Now,  if at least one of the two sub-processes tends to infinity with strictly positive probability, then {\it both}
must do so according to the definition of $\alpha$ ($0<\alpha<\infty$)
as the a.s. limiting ratio of $\G_t^i/\G_t^h$ as $t \to \infty.$ Hence, recalling (17) of part (i), it suffices to show that 
$$P\Big(\G_t^h \to \infty \Big|m_hF_h(\tau)>1\Big)=1-P\Big(\G_t^h \to 0\Big|m_hF_h(\tau)>1\Big)>0.$$ Since, by Definition 6.1, the process $(\G_t^h)$ is a RDBP this follows however already from
Theorem 4.4 ii) b) of Bruss and Duerinckx  (2015).

\smallskip
This completes the proof of Theorem 6.1.
\qed

\bigskip

\noindent{\bf 6.2 Remarks}

\smallskip
\noindent 1. Note that in (i) we neither need nor prove that the value $\alpha$ is unique but only that there is a one-to-one correspondence between $\tau$ and $\alpha,$ if such a candidate value $\alpha$ exists. Indeed, in Subsection 6.5.2 we will give examples with several candidates $\alpha$ for an equilibrium. 

\smallskip\noindent
2. In the proof of part (i) and at the end of part (iii) we did not use the Markov property of $(\G^h_t)$ and $(\G^i_t)$ although this would have shortened the proof. Indeed, we do have the Markov property  from the assumption that both processes are RDBPs (see Prop. 4.1 of Bruss and Duerinckx (2015)). However, we only used that $n_h$ ($n_i$) individuals in the home-(immigrant)-population will have no descendants with probability at least $p_0^{n_h}$ ($p_0^{n_i}$). The statement (17) holds more generally and may leave room for introducing more general processes but this direction is  not pursued in the present paper.

\subsection{The role of $\tau$ as a threshold claim}

It is the value $\tau$ which may be seen as an approximate upper threshold claim for any individual in the combined process $\G(t),$ provided the effectives of the sub-populations are not too small. Individuals claiming less than $\tau$ will have a good chance to remain in the society and to reproduce whereas those claiming more than $\tau$ will not. In any finite state at time $t$ the true threshold will be the value $\tau_t:=\tau(\omega)$ solving, for the same $\omega$ equation (13).
The true value $\tau_t$ depends of course on the empirical distribution function of the merged list of claims submitted by both sub-populations. 

\smallskip
In the following Lemma we give a strong bound on the speed of convergence of $\tau_t \to \tau.$ 

\bigskip
\noindent {\bf Lemma 6.1} If the constraints (15) are satisfied with $m_hF_h(\tau)=m_hF_h(\tau) >1$ (strict inequality), then, conditioned on survival of both sub-processes, the threshold values  $\tau_t, \,t=1, 2, \cdots ,$ defined in (18) converge exponentially quickly to their limit $\tau$ figuring in equation (15).

\medskip
\noindent {\bf Proof.} The Dvoretzky–Kiefer–Wolfowitz inequality implies that the empirical distribution function $F_n$ for $n$
i.i.d random variables following the distribution function $F$ satisfies 
$$P(\sup_{x\in \R} |F_n(x)-F(x)| >\delta ) \le 2e^{-2 n \delta^2}$$ for all $\delta>0$. For $F:= F_h$, say, this implies in particular
$$\forall \delta>0:P\Big( \big|F_{h, D_t^h(\G_t^h)}(x)-F_h(x)\big|>\delta\Big|D_t^h(\G_t^h)\Big)\le 2e^{-2 D_t^h(\G_t^h) \delta^2} ~{\rm a.s.},$$ not depending on $x,$ and the corresponding inequality holds if we replace $F_h, D^h_t, $ and $\G^h_t$ by $F_i, D_t^i$ and $\G^i_t$ respectively.  Now, the distance between adjacent claims on the merged list of increasing order statistics of the claims at time $t$ cannot be larger than the distance between adjacent claims in any of the separate lists. Thus, if we denote the true and the empirical distribution function of the merged list of claims at time $t$ by $J_t(\cdot)$ and $\tilde J_t(\cdot),$ respectively, we obtain
$$\forall x>0, \forall \delta>0:P\Big( \big|\tilde J_t(x)-J_t(x)\big|>\delta\Big|\G_t^h, \G_t^i\Big)\le 2  e^{-2 (D_t^h(\G_t^h)\vee D_t^i(\G_t^i)) \delta^2} {\rm ~a.s.}$$ independently of $x,$ where $a\vee b$ denotes the maximum of $a$ and $b.$ 

Now, with both long-term multipliers $m_hF_h(\tau)$ and $m_hF_h(\tau)$ being strictly greater than one, conditioned on survival of both processes $(\G_t^h)$ and $(\G_t^i),$ the random variables $\G_t^h$ and $\G_t^h$ will tend exponentially quickly to infinity as $t$ tends to infinity.
 The same must then hold for the respective numbers of descendants $D_t^h(\G_t^h)$ and $D^i(\G_t^i).$ It follows from the preceding inequality that the functions $\tilde J_t(x)$ and $J_t(x)$ converge, as $t \to \infty,$ exponentially quickly to each other for all $x.$ 
 Since $J_t$ is a convex mixture of the two absolute continuous functions $F_h$ and $F_i,$ $J_t$ is itself absolutely continuous. Thus $J_t$ is also uniformly continuous on any compact interval containing $\tau=\lim_{t \to \infty} \tau_t$ as an interior point. Consequently, for any $\epsilon>0$ there exists a constant $L>0$ not depending on $x$ such that
$$\forall \epsilon >0: |\tilde J_t(\tau_t )- J_t(\tau)| < \epsilon \,\implies |\tau_t-\tau|< L\, \epsilon.$$ Hence the exponential speed of convergence of $\tilde J_t(\tau_t )\to J_t(\tau)~{\rm a.s.}$ carries over to the speed of convergence of $(\tau_t)\to \tau$ a.s.,  and Lemma 6.1 is proved.\qed

\medskip

\subsection{Interplay of natality, productivity and resource claims}
\medskip
One possible solution  $\tau$ of equation (16) can  directly be read off, namely if the same $\tau$ solves simultaneously the two equations
\begin{align} m_h \int_0^\tau xdF_h(x) =r_h; ~~m_i\int_0^\tau xdF_i(x) =r_i. \end{align} This solution may be seen as a case of {\it perfect self-sufficiency} of the two cohabitating sub-populations under the {\it wf}-policy of attributing resources. Indeed, (22) says in words that an individual produces in expectation exactly what its expected number of children with {\it accepted} claims will consume in expectation.

 However, this coincidence can hardly be hoped for in reality. For instance,  if  $F_h=F_i$ then we  have to assume $r_h/m_h=r_i/m_i.$ Immigration goes mostly from the poorer nation into the direction of the richer one so that, as richness is positively correlated with productivity, we typically have $r_i<r_h.$ Moreover, since natality rates are usually higher in poorer countries the case $m_i>m_h$ is again more typical. Thus the ratio $r_h/m_h$ is usually  larger, and often substantially larger, than the ratio $r_i/m_i.$ 
We conclude that for Model I the existence of an equilibrium is an exception rather than the rule. 

\medskip\begin{quote}
The
illustrations in Figure 1 and Figure 2 display in a simplified form the typical phenomenon in the case of beta-densities of the variables claims.\end{quote}
 \medskip
\begin{figure}[ht]
\centering
 \includegraphics[width=0.35\textwidth, angle=0]{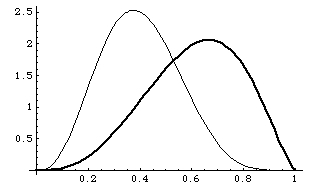}
\caption{Two beta-densities $f_h$ and $f_i$ on $[0,1]$ exemplifying in simplified form the distribution of claims from individuals of the home-population and immigrants: 
$f_h(x)$ with parameters $(a,b)=(4, 2.5)$ (in bold), 
and $f_i(x)$ with parameters $(a,b)=(4, 6).$}
\end{figure}

\begin{figure}[ht]
\centering
 \includegraphics[width=0.35\textwidth, angle=0]{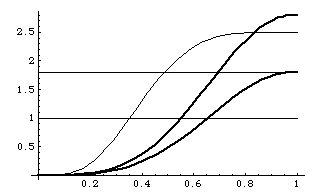}
\caption{This graph shows, for $m_h=1.8$ and $m_i=2.5,$ the functions $m_h\,F_h$ and $(m_h+1)\,F_h$ (both in bold), as well as $m_i\,F_i$ (in grey).}\end{figure} 
\begin{quote}The graphs of $m_h\,F_h$  and $m_i\,F_i$ with $m_h<m_i$ do not intersect on $[0,1],$ and hence the constraints (15) cannot be satisfied.  This exemplifies the  situation of immigrants having  higher birth rates and lower productivity and thus lower claims than individuals in the home-population. Possible control actions to allow for an intersection above the level 1 would be to try to increase $m_h$ (upper curve)\end{quote}

\medskip
\noindent Control actions of the type of increasing the natality of one sub-population or reducing the natality of the other one, as indicated in Figure 2, are mathematically easy to understand, but they are arguably difficult to accomplish in practice because, often enough, family size standards are strongly anchored in the culture or tradition of sub-populations.

\subsubsection{Self-sufficiency}
It is customary to refer to a population as being {\it self-sufficient} if, with respect to natality, productivity, and consumption,  it would be able to live on its own. Let us make this more precise by requiring that, if all parameters stayed the same forever, it would  survive on its own forever with positive probability.
What then would self-sufficiency mean for our model (Model I), where both sub-populations live in co-habitation profiting from a common resource space? Can they converge to an equilibrium?

To answer this question we look again at Theorem 6.1. 
Equation (16)  shows that if the mean production $r_h$ of the home-population and $r_i$ for the immigrants are kept constant, then the l.h.s. is increasing in $\tau$ for any pair of distribution functions $F_h$ and $F_i.$ Then, as we have seen, for each
$\alpha$ there can exist at most one $\tau:=\tau(\alpha),$ and vice versa. Solving equation (16) for $\al$ yields \begin{align} \al=\al(\tau)=\frac{r_h-m_h\,\int_0^\tau x dF_h(x)}{m_i\,\int_0^\tau x dF_i(x) -r_i},\end{align} where
$\alpha$ should satisfy $0<\alpha<\infty.$ 
Equation (23) has a non-trivial practical aspect worth being stated in words, with a trivial proof.  Let us recall for this the notion of perfect self-sufficiency we introduced at the beginning of subsection 6.4 for the case that a single value $\tau$ solves simultaneous both equations in (22). Then we have

\medskip
\noindent
{\bf Corollary 6.1}: Unless in the case of perfect self-sufficiency  an equilibrium can only exist if one sub-population becomes a net contributor to the Society's resource space, and the other one a net consumer of it.

\medskip\noindent
{\bf Proof}:  By  definition of an equilibrium we must have $\al \, \in \, ]0, \infty[,$ and hence  the numerator and denominator in (23) must have the same sign. Hence, if the sign is well-defined (meaning different from $0$ in most countries), the statement is true. If  both numerator and denominator vanish in equation (23), this implies that the same $\tau$ solves simultaneously both equations in (22).  \qed.

\medskip \noindent {\bf Remark 6.3} As far as one can judge from the media, Corollary 6.1, as trivial as it is, seemingly contradicts the intuition of contemporary decision makers. The wide-spread feeling is that if sub-populations
are "doing sufficiently well" in the sense that they could each live on their own, everything will be fine for a peaceful cohabitation. This is not true on the level of a long-term equilibrium. Without further control, the latter can only be attained for {\it  perfectly self-sufficient} sub-populations  in the sense of (22). Otherwise, without control, one sub-population will take over. The point is that there is little reason to believe that an adequate control would establish itself.

\subsubsection{Lorenz curves}

Self-sufficiency of a population does not only depend on its total resources but, to some extent,  also on the distribution of resources among the population. For instance, if the poorest class of individuals owns a relatively even smaller part of the total wealth it may be inclined to leave and, to close the gap, have to be replaced by more expensive individuals. 

The distribution of wealth is usually presented in a Lorenz curve. Accepted claims of individuals in our model can be seen as  increments in the Lorenz curve, and  Jacquemain (2017) showed that certain results of Bruss and Duerinckx (2015) allow indeed for  a clear re-interpretations in these terms.  Economists  may be interested in finding a link between (23) and the two Lorenz curves for the sub-populations.  Although the mean family sizes $m_h$ and $m_i$ intervene also in (23), this link remains worth studying. We should also mention that
Wajnberg (2014) interprets those results more freely, but some of these interpretations seem harder to justify.

\subsection{Multiple possible equilibria}
So far we have seen that, if an $\alpha$-equilibrium exists, then the solution $\tau:=\tau (\alpha)$ solving 
equation (15) is unique. When we introduced $\alpha$ in Definition 6.2 we introduced it as a random variable, implying that we may have several candidates for an equilibrium.  When will this be the case? 

We first look at a simple model for comparison.
\subsubsection{Comparing Model I with two Galton-Watson processes in co-habitation} Consider two Galton-Watson processes  $(Z_t^{(1)})_{t=1,1, \cdots}$ and $(Z_t^{(2)})_{t=1,2, \cdots}$ with respective  reproduction means $m_1$ and $m_2$.
We assume the usual conditions $p_0^{(j)}>0,\, p_0^{(j)}+p_1^{(j)}<1$ for $j=1, 2$. Recall also that we supposed for all random variables in this paper the existence of second moments, so that in particular the inequalities $\E(Z_1^{(j)} \log{Z_1^{(j)}}|Z_0^{(j)}=1)< \infty$ hold for $j=1,2$. It is well known that in this case the processes $(Y_t^{(1)})$ and $( Y_t^{(1)})$ defined by $$Y_t^{(1)}=Z_t^{(1)}/m_1^t,~~Y_t^{(2)}=Z_t^{(2)}/m_2^t, ~t=1, 2, \cdots$$ are a.s.-converging martingales so that
$Y_t^{(1)}/Y_t^{(2)}$ converges a.s. to a random variable. If $m_1\not =m_2$ then clearly only a degenerate limit $0$ or $\infty$ can exist for $(Z_t^{(1)}/Z_t^{(2)})_{t=1, 2, \cdots}$. If $m_1=m_2$ however, then the process $(Y_t^{(1)}/Y_t^{(2)})$ coincides with
$(Z_t^{(1)}/Z_t^{(2)})$ and thus converges a.s. to an equilibrium in the sense of our Definition 6.2 where this limit is distributed like the ratio of two functions of normal distributions (see e.g. Hall and Heyde (1980), subsection 1.3). Now, once both sub-populations $(Z_t^{(1)})$ and $(Z_t^{(2)})$ are sufficiently large, it follows from the i.i.d. reproduction of individuals that 
\begin{align}\frac{Z^{(1)}_{t+k}}{Z^{(2)}_{t+k}}
\sim \frac{Z^{(1)}_{t} m_1^k}{Z^{(2)}_{t} m_1^k}
=\frac{Z^{(1)}_{t}}{Z^{(2)}_{t}}, ~{\rm~ as}~ k \to \infty\end{align} 
and from the strong law of large numbers that the conditional distribution of $Z_{t+k}^{(1)}/Z_{t+k}^{(2)},~ k=1, 2, ...~$ given $Z_{t}^{(1)}$ and $Z_{t}^{(2)}$ concentrates around $ Z_{t}^{(1)}/Z_{t}^{(2)}.$ In simplified language we may say that the earlier history of states of the two Galton-Watson processes points  quickly to some sufficiently small neighbourhood of the equilibrium $\alpha_Z=\lim_{t\to\infty}Z_t^{(1)}/Z_t^{(2)}$ for which there are uncountably many candidates.

\subsubsection{Example of multiple candidates for equilibria}
Returning to Model I we see the similiarity with the long-term reproduction rates $m_hF_h(\tau)$ and $m_iF_i(\tau)$ which have to coincide in order to allow for an equilibrium. At the same time we also see that  RDBPs are  much more restrictive on possible equilibria because we now need a value
$\tau$   satisfying $m_hF_h(\tau)=miF_i(\tau)\ge1$, and, at the same time for a candidate $\alpha,$  the corresponding solution of equation (15). In the interesting case where $F_i$ and $F_{h}$ do not coincide on an interval of positive Lebesgue measure, where not both are equal to zero or one at the same time, it is clear that, in contrast to the model of two independent GWPs, there cannot exist intervals of possible equilibria but at most countably many candidates. 

Now, intuitively, we should be able to find several candidates if we choose $F_h$ and $F_i$ as well as $m_h$ and $m_i$ sufficiently close to each other with several intersection points of $m_hF_h=m_iF_i$ on the set of $t$ such that $m_hF_h(t)\ge1$ and then choose parameters $r_h$ and $r_i$ to keep this compatible with equation (15). Indeed, following this intuition leads easily
to an example showing multiple possible equilibria in Model I. Here is one:

\medskip \noindent Let $m_h=3$ and $m_i=2.8,$ and let for $n\in \{1,2,\cdots\},$$$F^{(n)}_h(t)=t+(2n\pi)^{-1}{\rm sin}(2n\pi t); ~F^{(n)}_i(t)
=t-(2n\pi)^{-1}{\rm sin}(2n\pi t).$$ For  $n=4$, for example, we obtain for $r_h=1$ and $r_i=0.5$ three candidates which  (using {\it Mathematica}, rounded to four decimals) are the points $(\tau_1=0.6031, \alpha(\tau_1)=13.3675), ~(\tau_2=0.7795, \alpha(\tau_2)=0.3681)$, and $ (\tau_3=0.8424, \alpha(\tau_3)=0.0016).$ Similarly as in our argument shown through (24) we would see, as soon as we can observe the processes $(\G_t^h)$ and $(\G_t^i),$ which of these equilibrium candidates is relevant. Note how different these values can be.\subsubsection{Importance in practice}
The preceding example is artificial in the sense that $F_h$ and $F_i$ mimic each other with one being a very similar delayed version of the other one. 
For larger $n$ we may find more (isolated) candidates. 
 However, even if such examples made sense in a real-world problem, this would hardly lead to a confusion  because the two sub-populations would be typically 
in the thousands or millions. With the proven exponential speed of convergence (see Sub-section 6.3), the current  states will determine with probability close to one the relevant $\tau,$ and thus the relevant equilibrium. 
Moreover, it is hard to imagine realistic situations where $m_hF_h$ and $m_iF_i$ allow for several points of intersection. The phenomenon of multiple equilibrium is of mathematical interest rather than relevant in practice, and the typical real world problem coming with immigration is to find an accessible equilibrium rather than having several  candidates. 

\subsubsection{Open problem, and its impact.}
The reader will have noticed that Theorem 6.1 would be a necessary and sufficient criterion for an equilibrium to exist with strictly positive except that we have not shown that the slightly stricter constraints $m_hF_h(\tau)>1$  and $m_iF_i(\tau)>1$ are also necessary. Is it thinkable that $m_hF_h(\tau_t)$ and $m_iF_i(\tau_t)$ converge,  as $t\to \infty,$ so slowly to $1$ that both sub-processes may still tend to infinity and thus enable the existence of an equilibrium?  In the case of limiting criticality the $\tau_t$ seem difficult to control, and the author must leave this question open. 

Independently of this, it is the necessary conditions which should attract our main interest for applications since it is the necessary steps for survival which must be taken (ad-hoc) by the decision makers in each generation. What would it mean to say that "if the current conditions of the RDBPs would stay forever, then such or such condition would also be sufficient for the existence of an equilibrium"? The necessary conditions must be observed by the society obligation principle, and unpredictable events leave little motivation to study  sufficient conditions for survival in a long-term probability model which one cannot predict, not even approximately. The mathematical side of the open question is a  challenge (also encountered in Bruss and Duerinckx (2015)), but concerning applications, the author sees no impact.

\section {The effect of integration of immigrants}
As we can see from Theorem 6.1, the constraint qualification  $m_h\, F_h(\tau)=m_i\, F_i(\tau) \ge1$ is demanding, and is likely to be a major obstacle for the existence of a real-world equilibrium. As argued before, immigration goes  usually (apart from politically persecuted individuals) from poorer countries with higher birth rates, i.e. $m_i>m_h$, into a richer one with stochastically larger claims. 

What can be hoped for from integration, that is,  immigrants adapt sufficiently quickly the economically relevant behaviour of the home-population? By {\it economically relevant} we refer
to having the same parameters of mean natality $m_h$ and mean productivity $r_h$, and the same distribution function $F_h$ of claims, whereas cultural
or religious factors are not taken into consideration.
Different assumptions about the mechanisms of integration lead then to different models. We study a model built on, what we call, {\it fractional integration.}

\subsection{Fractional Integration - Model II}
 Suppose that in each generation $t$, the fraction $\vfi$ of the number of those individuals currently seen as belonging to the immigrant-population,  will integrate into the home-population in the sense that, once integrated, they share the same reproduction parameter $m_h,$ the same mean productivity $r_h,$ and the same claim distribution function $F_h.$ The corresponding bivariate process will be denoted by\begin{align}\left(\Gamma^{(\vfi)}(t)\right)_{t=1, 2, \cdots}=\left(\Gamma_t^{h,(\vfi),}, \Gamma_t^{i,(\vfi)}\right)_{t=1, 2, \cdots}.\end{align} We note that, from a formal point of view, we should append the fraction $\vfi$ as an additional index in $(\G_t)$, $(\G_t^h)$ and $(\G^i_t)$, because these are now also functions of $\vfi.$
 
\medskip \noindent
{\bf Definition 7.1}: The bi-variate process $\left(\Gamma^{(\vfi)}(t)\right)_{t=1, 2, \cdots}$consisting of the home-population and the immigrant-population when, in each generation $t$ a fraction $\vfi$ of the immigrant-population integrates into the home-population will be called {\it $\vfi$-fractional integrated} process. We refer to the process $\left(\Gamma^{(\vfi)}(t)\right)_{t=1,2, \cdots}$ with constant integration fraction $\vfi$  as Model II. 

\medskip\noindent
In the spirit of Definition 6.1 we  define:

\medskip
\noindent{\bf Definition 7.2} We say that the bivariate process $(\G^{(\vfi)(t)})_{t=1,2, \cdots}$ tends to an equilibrium if there exist a value $\vfi \in [0,1[$ and a corresponding random variable $0< \alpha_\vfi < \infty$ such that
\begin{align}~P\left( \lim_{t \to \infty} \frac{\G_t^{i, \,(\vfi)}}{\G_t^{h, \,(\vfi)}}  = \alpha_\vfi\Big|\G_t^{i,(\vfi)}\not\to 0, \G_t^{h,(\vfi)}\not\to 0\right)=1.\end{align} 

\smallskip

\noindent As we shall see in the following Subsection it leads in general to no ambiguity to think of $\vfi$ as being fixed and, for simplicity of notation,  to drop the upper index $(\vfi)$ in the
sub-processes whenever $\alpha_\vfi$ and $\vfi$ are implicit functions of each other with a unique solution, and this will always be the case if the limit $\alpha_\vfi$ exists.
In particular, $\alpha_0=\alpha,~\G_t^{i, \,(0)}=\G_t^i$ and $\G_t^{h, \,(0)}=\G_t^h. $ To save space in the longer equations to come, we define 
$$\Phi_h(t)=\int_0^t x dF_h(x);~~\Phi_i(t)=\int_0^t x dF_i(x).$$
Our objective is now to derive the corresponding equilibrium conditions. 

\subsection{Equilibrium equations for the combined process with fractional integration}
The random total resource space at time $t+1$ becomes then in the simplified notation described above
\begin{align}
\tilde R\Big(\G_t^h+\vfi \G_t^i, (1-\vfi)\G_t^i)\Big)=R^h(\G_t^h+\vfi \G_t^i)+R^i((1-\vfi)\G_t^i).
\end{align}
If $\vfi \G_t^i$ is integer-valued, then $R^h(\G_t^h+\vfi \G_t^i)$ is well-defined and follows the same distribution as $R^h(\G_t^h)+R^h(\vfi \G_t^i).$ To exclude ambiguity we think of $\vfi \G_t^i$ as being defined as $ [
\vfi \G_t^i]$, where $[x]$ denotes the floor of $x$, say. This is asymptotically of no importance if $\alpha _\vfi$ exists, of course, and will therefore no longer be mentioned.

\medskip

\noindent{\bf Lemma 7.1}
The random BRS-equation for the $\vfi$-integration process equilibrium is given by
\begin{align} \Big( D_t^h(\G_t^h)+D_t^h(\vfi \G^i_t)\Big)\Phi_h(\tau_t) + \Big(D_t^i(\G_t^i)-D_t^i(\vfi\,\G_t^i)\Big)\Phi_i(\tau_t)\\
= R^h\Big(\G_t^{h}\Big)+R^h\Big(\vfi\G_t^{i}\Big)+R^i\Big((1-\vfi)\G_t^{i}\Big),~~t=1, 2, \cdots\end{align} and there can exist  at most one solution $\tau_t$.

\bigskip\noindent
{\bf Proof:} We have first to show that the equation defined by (28)=(29) is the random BRS-equation for Model II. We see that it is well-defined since it is well-defined for all $\omega\in \Omega$  and all $\vfi\in [0,1].$

To understand the terms in (28) and (29), recall that there is a shift of the fraction $\vfi$ of the immigrant-population into the home-population  from generation $t$ to $t+1.$ By our assumption these individuals now reproduce and  consume (claim) independently like the other individuals of the home-population, where $\Phi_h(\tau_t)$ is the corresponding average accepted claim. This yields the first product in (28), which corresponds to the  total amount of accepted claims of resources submitted by the current home-population.  

The second product in  (28) reflects the corresponding reduction of the number of individuals, and thus of their total claim for the immigrant-population. 

The r.h.s. of the equation, that is  (29), gives accordingly  the  balance of the random total resource space contributions. Here we have used throughout the additivity of resource production, the i.i.d. reproduction within the same sub-population,  and the definition of $\vfi$-integration. Hence, according to equation (12), this is the random BRS-equation of Model II.

\smallskip
To see that for given $\vfi, \alpha_\vfi$ there exists at most one solution $\tau_t,$ recall that $F_h(x)$ and $F_i(x)$ are absolute continuous functions so that
 $\Phi_h(x)$ and $\Phi_i(x)$
are also absolute continuous. Moreover, these four functions are strictly increasing in $x,$ so that we can repeat the arguments given in part (ii) of the proof of Theorem 6.1.  Therefore, for fixed natality  and production  laws  governing $D_t^h, D_t^i,  R_t^h$ and $R_t^i$,  there exists at most one solution $\tau_t,$
depending on, and well defined, for each $\omega\in \Omega$.
\qed

\bigskip \noindent
{\bf Theorem 7.1} For an equilibrium in Model II with integration rate $\vfi$ to exist it is necessary that there exists values $\tau>0$ and $\alpha_\vfi$ with $0<\alpha_\vfi<\infty$ satisfying the equation
\begin{align} m_h\Big( 1+\vfi \alpha_\vfi \Big)\Phi_h(\tau)+m_i \alpha_\vfi  \Big(1-\vfi\Big) \Phi_i(\tau)=r_h+r_i\alpha_\vfi + \vfi \alpha_\vfi
\Big(r_h-r_i\Big)
\end{align}
subject to the constraints \begin{align}m_h(1+\alpha_\vfi \vfi)F_h(\tau)=m_i (1- \vfi)F_i(\tau) \ge 1.\end{align}

\medskip\noindent{\bf Proof:} The proof is based on the facts that we can again apply the strong law of large numbers within each sub-population, and on the hypothesis that the limit $\alpha_\vfi=\lim_{t \to \infty} \G_t^{i, (\vfi)}/ \G_t^{h, (\vfi)}$ exists. 
\smallskip
The proof is therefore similar in its structure as the proof of Theorem 6.1., and we can be more concise. 

We show for example the limiting result on the r.h.s. of the equation, that is (29). Dividing both sides of the equation in Lemma 7.1 by
$\G_t^{h}:=\G_t^{h, (\vfi)}$ we can rewrite the r.h.s. of the new equation in the form
$$\frac{R^h(\G^{h}_t)}{\G_t^{h}}+\frac{R^{h}(\vfi \G_t^{i})}{\vfi \G^{i}_t}~\frac{\vfi \G_t^{i}}{\G_t^{h}}+\frac{R^{i}((1-\vfi)\G^{i}_t)}{(1-\vfi)\G^{i}_t}~\frac{(1-\vfi)\G^{i}_t}{\G^{h}_t},$$ which, according to our assumption of independence within each sub-population, converges a.s. to $r_h+r_h \vfi \alpha_\vfi+r_i(1-\vfi)\alpha_\vfi$, as $t$ tends to infinity. The latter is the r.h.s. of (30). 

Since the r.h.s. allows for a limit conditioned on survival, the limit on the l.h.s. must also exist and, of course, coincide.  Hence in particular, under the same condition of survival, we must have $\tau_t\to \tau$ for some  $\tau.$  It is straightforward to check that  the l.h.s. of (28) divided by $\G_t^t$ yields then the l.h.s. of the limiting equation (30).

\smallskip
The proof of the constraint qualifications in (31)  is also quite similar. Let us show this for the home-population.

Recall that the probability of an individual in the home-population having no offspring equals $p_0^h>0.$ Consequently, given $\G_t^{h,(\vfi)}=n,$ the absorbing state $0$ is accessible within the home-population with at least probability $(p_0^h)^n>0$. As we have seen by the Borel-Cantelli Lemma type argument given in the proof of part (iii) of Theorem 6.1, the sequence of conditional expectations $$\left(\E(\G_{t+1}^{h,(\vfi)}|\G_{t}^{h,(\vfi)}>0)\right)_{t=1,2, \cdots}$$ must therefore not stay bounded  because otherwise $(\G_t^{h,(\vfi)}) \to 0$ almost surely as $t \to \infty. $ It follows that the number of descendents in generation $t+1$ on the home-population must tend to infinity as $t$ tends to infinity and thus behave asymptotically like
$$m_h\left(\G_t^{h,(\vfi)}+\vfi \G_t^{i,(\vfi)}\right), \,t=1, 2, \cdots.$$
Since the probability of a claim of an individual in the home-population to be accepted equals $F_h(\tau_t)$
their total number in generation ${t+1}$ behaves like
$$\G_{t+1}^h\sim m_h\left(\G_t^{h,(\vfi)}+\vfi \G_t^{i,(\vfi)}\right) F_h(\tau_t), \,t=1, 2, \cdots.$$ 
The arguments for the immigrant-population follow the same line of reasoning, except that the immigrant-population looses the fraction  $\vfi$ of its current effectives to the home-population.
Its number of descendants in generation $t+1$  behaves thus asymptotically like $m_i\left(\G_t^{i,(\vfi)}-\vfi \G_t^{i,(\vfi)}\right)$, so that
$$\G_{t+1}^i\sim m_i\left(\G_t^{i,(\vfi)}-\vfi \G_t^{i,(\vfi)}\right) F_i(\tau_t), \,t=1, 2, \cdots$$Hence the necessary conditions for the existence of an equilibrium seen in (15) hold for the $\vfi$-integrated process correspondingly 
in the form 
\begin{align}m_h(1+\alpha_\vfi \vfi)F_h(\tau)\ge 1 \rm{~and~} m_i (1- \vfi)F_i(\tau) \ge 1.
\end{align} (We note here a slight asymmetry in the sence that $\alpha_\vfi$ does not appear in the second inequality in (32)).
 Finally, going then through the steps (19) to (21) in an analogous way with the new factors shows that  both factors,  $m_h(1+\alpha_\vfi \vfi)F_h(\tau)$ and $m_i (1- \vfi)F_i(\tau),$ which are the asymptotic multiplication factors for the home-population, respectively,  for the immigrant-population, must again coincide in order to allow for an $\alpha_\vfi$-equilibrium.

This will complete the proof.\qed

\bigskip
\noindent{\bf Remark 7.1} If we replace both inequalities in (31) by "$>1$" then, conditioned on the event $\{\G_t^{i,(\vfi)}/\G_t^{h,(\vfi)}\not\to 0\}\cap \{\G_t^{i,(\vfi)}/\G_t^{h,(\vfi)} \not\to \infty\},$ Theorem 7.1 states a sufficient condition for an equilibrium to exist with positive probability. Using the new asymptotic multiplying factors, the proof follows exactly the same lines as (iv) of the proof of Theorem 6.1. \qed

\subsubsection*{Illustration of an equilibrium in Model 2}
\medskip Equation (30) in Theorem 7.1 is  linear both in $\vfi$ and in $\alpha.$  The choice of solving this equation for $\vfi$ or $\alpha$ will primarily depend on our main objective: Is our question, first of all, which equilibria are, in principle,  feasible with a certain upper bound for $\vfi$, or alternatively,  is the relevant question rather what integration fraction  $\vfi$ would be needed to obtain an equilibrium at a given desired level $\alpha$?

\smallskip
\noindent Suppose we decide for the latter, that is, we solve (30)
for $\alpha.$ This yields

\begin{align}\alpha:=\alpha(\vfi, \tau) = \frac{r_h-m_h\,\Phi_h(\tau)}{\vfi m_h \Phi_h(\tau)+(1-\vfi)m_i\Phi_i(\tau) - r_i +\vfi(r_i-r_h)}\end{align}

\medskip
This equation describes thus a surface over ${\cal S}_{\tau} \times {\cal S}_{\vfi}$, where ${\cal S}_{\tau}$ and ${\cal S}_{\vfi}$ denote the domains of $\tau$ and $\vfi$, respectively. The graph of $\alpha(\vfi, \tau)$ tends to plus or minus infinity along the curve described by all those points $(\vfi,\tau)$ in which the denominator of $\alpha(\vfi, \tau)$ vanishes. It seems therefore easier and sufficiently informative  to look only at the graph representing the constraints specified in (31). 

\smallskip
If we fix $\alpha,$ say, we obtain then three surfaces of the form $Z=Z(\vfi,\tau)$, of which one is the plane $Z_1 \equiv 1$ neither depending on $\vfi$ nor on  $\tau,$ and the others being $Z_2(\vfi,\tau)=m_h (1+\alpha \vfi)F_h(\tau)$ and  $Z_3(\vfi,\tau)=m_i(1-\vfi)F_i(\tau).$ 

\medskip
 We now look at Figure 3 based on our examples  of beta-distributions $F_h$ and $F_i$ presented in Figure 2. In Figure 3, $\alpha$ is fixed. We have chosen $\alpha=1/2$. Note that the corresponding surface  $\alpha(\vfi,\tau)$ cannot be plotted in a meaningful way in the same graph because the constraints surfaces depend on $\alpha.$ The intersection of $Z_2$ and $Z_3$ visible above $Z_1,$ i.e. the path going up from the plane to the upper left side, is thus the image of the points satisfying
 the constraints given in (31). The projection of this path on ${\cal S}_{\tau} \times {\cal S}_{\vfi} = ]0,1[^2$ intersected with the projection of the level-curve $\alpha=1/2$ (both not visible here) present then the set of $(\vfi,\tau)$ allowing for an $\alpha$-equilibrium.
 
For $\vfi=0.20$  for example,
numerical computation (Mathematica) shows that $\tau= 0.80106...$. 

\medskip
\begin{figure}[ht]
\centering
 \includegraphics[width=0.35\textwidth, angle=0]{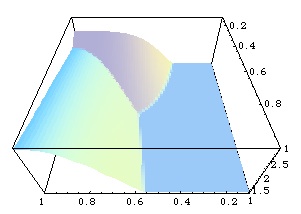}
\caption{Illustration of Constraints for $\alpha =1/2$ and $0<\vfi<1, 0 < \tau<1$.}\begin{quote}  In this graph we show for $\alpha=1/2$  the constraints as a function of $\vfi$ and $\tau$ for $F_h$ and $F_i$ presented in Figure 1 and Figure 2. The other parameters are $m_h=1.8, m_i=2.5, r_h=1.0, r_i=0.6.$  The graph shows  the plane $z \equiv1$ (in blue) and the two surfaces corresponding to the two constraints. For a better view, the graph has been rotated to the view-point $(0,1,1)$. \end{quote}\end{figure}
\subsubsection{Conclusions for Model II} The introduction of integration makes a clear  and important difference compared with Model I without integration. The constraint qualifications (30) are, at least in principle, much easier to satisfy than those for Model I we saw in (15). Hence, if the home-population can afford the investments needed for a more generous integration rate $\vfi$ it will frequently succeed to obtain an equilibrium. Given that the mean reproduction  $m_h$ is,  in real world, essentially smaller that $m_i$ and that it is not easy,  to overcome traditional and cultural differences explaining larger differences of parameters (also with respect to $r_h$ and $r_i$) we see that integration is an efficient factor of control. Having said this we understand of course that, in the real-world setting, a larger integration rate $\vfi$ may change many characteristics of the home-population.

We also mention that the sensitivity of equilibrium points observed in Model I (see Sub-section 6.5) with respect to the parameters and/or claim distributions can also be observed in Model II with respect to the additional integration  parameter $\vfi.$ This is hardly a surprise because fractional integration can  be interpreted as a change of the effective reproduction of the two sub-populations.

\section{General equilibrium equations}
 Studying the effect of immigration in a realistic way means however more. 
 We must assess,  at the same time, the effect of different parameters and different consumption features of two sub-populations, the effect of integration of one sub-population into the other one, and, in addition and in particular, the effect of an ongoing stream of new arrivals into one sub-population.
 
\smallskip
To reach this goal we should complete the model by allowing, in one form or another, new immigrants in each generation. 

Formally, we now define a tri-variate process
\begin{align}\left(\G^{(\vfi,I)}(t)\right)_{t=1, 2, \cdots}=\left(\Gamma_t^{h,(\vfi,I)}, \Gamma_t^{i,(\vfi,I)}, I_t\right)_{t=1, 2, \cdots}\end{align}  where $I_t\ge 0$ denotes the number of new immigrants in generation $t.$  If $(I_t)\equiv 0$ and $\vfi,$ as before, the fraction of members of the 
immigrant-population which integrates into the 
home-population, then the process $\left(\G^{(\vfi,I)}(t)\right)$ is defined as the bi-variate process defined in (25). 

\smallskip
To be consistent in our notation we denote the mean reproduction and the mean resource productivity per {\it new} immigrant by $m_{ni}$,  
respectively $r_{ni}.$ The claim distribution function for new immigrants is denoted by $F_{ni}(x),$ and,  correspondingly, we put $$\Phi_{ni}(x):=\int_0^t x dF_{ni}(x). $$

The interaction of $(I_t)$ with the two sub-processes can be modelled in many ways, and each model may have its own justification.
 It would, for instance, be easy to consider new immigrants directly as an integral part of the immigrant population with the same parameters $m_{ni}=m_i$ and $ r_{ni}=r_i,$ and the same distribution function of claims $F_{ni}(x)=F_i(x)$. Indeed, in this case only the effectives in the second component would change, i.e. $\G_t^i$ would become 
 $\tilde \G_t^{(\vfi),i}:= \G_t^{(\vfi),i}+I_t.$  However, this would in general not be a convincing setting because new immigrants arriving somewhen in $]t,t+1[$ cannot be expected to contribute to the resource space as much as those born at time $t$. This is why we  should allow for more freedom in the model. 
 
 \subsection {Modelling aspects}
As we shall see in the following, the approach we have proposed is flexible. Even keeping all factors different will cause no problem for our approach as long as either $I_t/\G^h_t$ or $I_t/\G^i_t$, conditioned on survival of both processes $(\G_t^h)$ and $(\G_t^i),$ can  be supposed to tend almost surely to some limit. We call this in brief {\it ni-limit condition}.

The hypothesis of the existence of such a {\it ni}-limit may be  restrictive unless the limit zero is permitted, and, indeed, this is what we do permit!  The condition  becomes then rather mild. Also, and in particular, it is a  reasonable condition because at least one of the sub-populations will typically decide how many new immigrants will be allowed to enter. Clearly, we do not have to specify which one, because, if an equilibrium exists then, if one {\it ni}-limit exists, both will exist.

The advantage of our approach is that to each model  corresponds then
a unique tractable BRS-equation, a unique BRS-inequality, and under the {\it ni}-limit condition a unique necessary condition for the existence of an equilibrium. All steps leading to the equilibrium equation will essentially follow the main scheme. The reader may agree that this unified structure allowing to pass from the simplest model (Model I) to a model with integration (Model II), and then up to the comprehensive model allowing in each generation new immigrants (which we will call Model III) is transparent and adds to  the tractability of the approach. 

\subsection{Fundamental equilibrium equation}
Having said that we can think of  many different models of how new immigrants should be assessed, the author thinks that, as a first approximation,  the following model  appeals reasonably well to reality. We name it here a {\it fundamental model} since it is our first model to include all essential processes concerned by  immigration into a new environment. 

Of course, we are fully aware that some modifications of what we propose as Model III, might attract more interest.  Given the set-up of our approach, several modifications would be equally tractable.\begin{quote} 

{\bf Fundamental model} (Model III):

(i) In generation $t$, that is somewhere in $[t, t+1[$, $I_t$ new immigrants join the co-habiting home-population and immigrant population.  We suppose that the $I_t$ may depend on $\G_t^h$ and/or $\G_t^h$ and that the process $(I_t)$ satisfies the $ni$-limit condition where a $ni$-limit $0$ a.s. is permitted. 

(ii) New immigrants are allowed to be  different from immigrants of the second or a later generation, and also  different from individuals belonging to the home-population. In generation $t$, the $I_t$ new immigrants have the right to consume (and do so according to the law $F_{ni}$)
but their production of new resources may be practically 0 or even negative.
The same is supposed to hold for their descendants
in the residual time before time $t+1.$ (This simplification seems justified since the residual time is likely to contain more individuals who consume than produce.) 

(iii) Immigrants present already for one generation or longer are still different from individuals in the home-population up to a random time until complete integration when they will become an integral part of the home-population. Descendants of $I_{t-1}$ are supposed to become members of the immigrant population but to still have their new immigrants consumption and production behaviour. 

(iv) Descendants of the generation $t-2$ and earlier have either disappeared, or are in $C_i$ or else already integrated into $C_h.$ 
\end{quote}

\noindent In conclusion, we allow for each of the three classes of individuals $C_h, C_i$ and $C_{ni}$ with (in general) different parameters of natality, productivity and different claim distributions. The $ni$-limit denoted by $\ell_{ni}$ is supposed to be a value determined by directives given by the home-population and/or the (established) immigrant population, and defined by $$P(\lim_{t \to \infty}(I_t/\G_t^h)=\ell_{ni}|\G_t^h\to \infty) =1.$$ 
 The process of integration is thought of as being fractional with parameter $\vfi$ introduced in Section 7.1.  We note that fractional integration leads in our setting to a geometric distribution of the random time of integration after immigration. As before, no intermediate steps of partial integration are considered. 

\smallskip\noindent The following Definition  Theorem will give the corresponding fundamental equilibrium equation.

\bigskip
\noindent{\bf Definition 8.1} The tri-variate process $(\G(t)^{(\vfi, I)})_{t=1,2, \cdots}$ defined in (34) is said to converge to an equilibrium, if there exists a value $\vfi\in [0,1]$ and a corresponding random variable $\alpha:=\alpha_\vfi$, such that conditioned on the survival of both sub-processes $(\G_t^h)$ and  $(\G_t^i)$
\begin{align} \lim_{t \to \infty} \frac{\G_t^{i, \,(\vfi, I)}}{\G_t^{h, \,(\vfi, I)}}  = \alpha_{\vfi}~{\rm a.s.}\end{align}

\smallskip\noindent In the preceding definition of an equilibrium between home-population and immigrant population it is formally not yet necessary to refer to the $ni$-limit condition and to the value $\ell_{\it ni}$ of the $ni$-limit. However, the latter will automatically intervene in the following criterion which displays the necessary equilibrium conditions.
To increase the transparency of the limiting BRS-equation, we present it in scalar product notation of vectors.

\bigskip\noindent
{\bf Theorem 8.1}  Let $\vec\Phi(t)$ denote the line vector function $\left(\Phi_h(t), \Phi_i(t),\Phi_{ni}(t)\right)$ and $\vec R$ denote the line vector $\left(r_h, \,r_i,\, r_{ni}\right)$. A necessary condition for the existence of an $\alpha$-equilibrium in  the tri-variate process 
 $(\G(t)^{(\vfi, I})_{t=1,2, \cdots}$ with $ni$-limit $\ell_{ni}$ is the existence of values $0< \vfi <1$, and $\tau>0$ solving the equation \begin{align}
\vec\Phi(\tau) \cdot
        \begin{pmatrix}
          m_h(1+\vfi\,\alpha) \\         
          m_i\left\{(1-\vfi)\alpha+G_{ni}(\tau)\right\}\}\\\ell_{ni}
          \end{pmatrix}=\vec R \cdot  \begin{pmatrix}
          1+\vfi\,\alpha \\         
          (1-\vfi)\alpha  \\
        G_{ni}(\tau)
          \end{pmatrix}    
  \end{align}
  where
  \begin{align*}G_{ni}(t):=F_{ni}(t) m_{ni}\ell_{ni}\left(F_h(t)m_h(1+\vfi \alpha)\right)^{-1},\end{align*}
and where $ F^{h}, F^{i}, F^{ni}$ and the parameters $m_i, m_h, m_{ni}, \ell_{ni} $ satisfy for $\alpha,$ $\vfi$ and $\tau$ the constraints
  \begin{align}
F_h(\tau) m_h \Big(1+\vfi\alpha \Big)= F_i(\tau)\,\Big\{m_i(1-\vfi)+G_{ni}(\tau)\Big\}\ge \,1.                
  \end{align}
  
  \bigskip
  \noindent{\bf Proof}: To avoid repetitions of mathematical arguments, we shall confine our proof to those parts which are different from the proof of Theorem 7.1.
  
We first show that the  random BRS-equation corresponding to (12) of Theorem 5.2 becomes now
for $t=1,2, \cdots$\begin{align} D_t^h(\G_t^h +\vfi \G^i_t)\,\Phi_h(\tau_t) + \left[~D_t^i\left((1-\vfi)\,\G_t^i\right)+F_{ni}(\tau_t)D_t^{ni}(I_{t-1})~\right]\Phi_{i}(\tau_t) +I_t\Phi_{ni}(\tau_t)
\\=~R_t^h\left(\G_t^{h}+\vfi\G_t^i\right)+R_t^i\left((1-\vfi)\G_t^{i}\right)+R^{ni}_t\left(F_{ni}(\tau_{t})D_{t-1}^{ni}(I_{t-1})\right). ~~~~~~~\end{align}
To see this we first note that
the first product on the l.h.s.\,of this  equation in line (38)
 does not change  compared with the first product in (28). The reason is that there are no direct transitions from the $I_t$ new immigrants (belonging to the class $C_{ni}$ ) into the home population, that is, no direct transitions from the class $C_{ni}$ into the class $C_h.$ The total random  consumption is thus $D_t^h(\G_t^h +\vfi \G^i_t)\,\Phi_h(\tau_t),$ and the limiting reproduction rate within $C_h$ becomes correspondingly
\begin{align}\lim_{t\to \infty} \frac{\G^h_{t+1}}{\G^h_t}=m_hF_h(\tau)(1+\vfi\alpha).\end{align}

\smallskip  Now, the descendants of the $I_{t-1}$ new immigrants of the preceding generation become part of $C_i.$  
 According to the assumptions (i)-(iv) of Model III, the random number $\vfi \G_t^i$ will again integrate into the home-population before re-producing and thus affect the first term.
Looking at the second term in (38) we see threfore an important change for the immigrant-population $(\G_t^i)_{t=1, 2, \cdots}$. 
The random number of descendants of the remaining  fraction $(1-\vfi)$ stays, as before, in $C_i.$ Now, however, the (random) fraction $F_i(\tau_t)$ of the descendants of $I_{t-1}$ will  submit claims in generation $t$ and thus also be  part of $C_i$,
thus counting for $\G_{t+1}^i.$  Multiplying the sum in brackets $[~~]$ in (38) with the random average consumption $\Phi_i(\tau_t)$ constitutes the second term and represents the random total consumption within the class $C_i$. 

The consumption of the new $I_t$ follows, individually, the law $F_{ni}$ so that the total  random consumption of the new immigrants equals $I_t\Phi_{ni}(\tau_t),$ which is the third term in (38).

\smallskip
Line (39) follows now correspondingly, except that we had supposed in Model III that there is no direct
contribution of resources by the $I_t$ new immigrants arriving in $[t,t+1[$.
  
 \medskip
 
 It remains to prove that the constraint qualifications
 in (37) are necessary for the existence of an equilibrium. 
 
We first note that the new immigrants $I_t$  intervene in the transition from  $(\G_t^h, \G_t^i)$ to  $(\G_{t+1}^h, \G_{t+1}^i)$ only on the side of consumption of resources but that the descendants of new immigrants from generation $t-1$  do intervene for $(\G_t^i)$ because they will also submit their claims in generation $t.$
The fraction $F_{ni}(\tau_t)$ of these will add  to the immigrant population. Further, the fraction $F_i(\tau_t)$ of their descendants will stay in the population. Hence, dividing the effectives in the class $C_i$ at time $t+1$ by the corresponding number in generation $t$ we obtain
\begin{align} \frac{\G_{t+1}^i}{\G_t^i}=\frac{1}{\G_t^i}\left(D^i_t\left((1-\vfi)\G_t^i \right) + D^{ni}_{t-1}(I_{t-1}) F_{ni}(\tau_{t-1})\right)\,F_i(\tau_t)\end{align}
To obtain the limiting reproduction rate in $C_i$ we now use
\begin{quote}
a)  The limiting equilibrium $\alpha_{\vfi,ni}$ is supposed to exist. Consequently the sequence $\tau_t$ converges to some value $\tau,$ and by continuity
all functions $F_h( \tau_t), F_i( \tau_t),
F_{ni}(\tau_t)$ and $\Phi_h(\tau_t), \Phi_i(\tau_t),
\Phi_{ni}(\tau_t)$ converge to their corresponding limits.

\medskip
b) The existence of~  $\ell_{ni}=\lim_{t\to \infty}(I_t/\G_t^h)>0$ implies the existence of the limit $\lim_{t\to \infty}
(I_{t-1}/\G^i_t)$ since, putting $\alpha=\alpha_{\vfi,ni}$ we have
$$\frac{I_{t-1}}{\G^i_{t-1}}\frac{\G^i_{t-1}}{\G^i_t}
\sim \frac{I_{t-1}}{\alpha \G^h_{t-1}}\frac{\G^i_{t-1}}{\G^i_t} \sim \frac{\ell_{ni}}{\alpha}\frac{\G^h_{t-1}}{\G_t^h}\sim\frac{\ell_{ni}}{\alpha}\left(m_h F_h(\tau)(1+\vfi\,\alpha) \right)^{-1},
$$ where the last step holds since the limiting multiplying factor in the class $C_h$ was already seen in (40) to be $m_hF_h(\tau)(1+\vfi\,\alpha).$
\end{quote}

\noindent Using a) and b), and again the strong law of large numbers, it is straightforward to show that the 
limiting multiplier of (41) becomes
\begin{align}
\lim_{t\to \infty}\frac{\G_{t+1}^i}{\G_t^i}=F_i(\tau)\left(m_i(1-\vfi) + \frac{F_{ni}(\tau) m_{ni}\ell_{ni}}{F_h(\tau)m_h(1+\vfi \alpha)} \right),
\end{align}
which proves the constraint qualification.

Note also that in our recursive setting, the number $I_{t-2}$ or the numbers of new immigrants from earlier generations have already been taken into consideration into the numbers $\G_{t-1}^i, \G_{t-2}^i \cdots$ and do not show up any longer.

 \smallskip To understand the additional term $G_{ni}(t)$  figuring in the equation (36) let us compute the limiting influence of the stream of new immigrants $(I_t)_{t=1, 2,\cdots}$ on both sub-populations. If we know it for one sub-population we then know it for the other one if we suppose that the limit  $\alpha_{\vfi, ni}$ defined in (35) exists. Keeping the delay of one generation for this influence in mind we
divide $I_{t-1}$ by $\G_t^h,$ (say) and obtain under the condition
 $\ell_{ni}=\lim_{t\to \infty} I_t/\G_t^h,$ 
\begin{align}\frac{I_{t-1}}{\G_t^h}=\left(\frac{I_{t-1}}{\G_{t-1}^h}\right)\left( \frac{\G_{t-1}^h}{\G_{t}^h}\right)\sim \ell_{ni}\, \lim_{t\to \infty} \frac{\G_{t-1}^h}{F_h(\tau_{t})\,D_{t-1}^h\left( \G_{t-1}^h+\vfi \G_{t-1}^i\right)}\\=
 \ell_{ni} \, \frac{1}{F_h(\tau) m_h (1+\vfi \alpha)},~~~~~~~~~~~~~~~~~~\end{align}
 Similarly, as in the proofs of Theorem 6.1 and Theorem 7.1, we can see that the limiting reproduction rate must be the same in both sub-populations in order to allow for an equilibrium, and our Borel-Cantelli type arguments apply again directly to see that that it must be greater or equal to one, completing the proof of Theorem 8.1. \qed
 
 \medskip
 \noindent {\bf Remark 8.1} If $\ell_{ni}=0,$ then it is intuitive that $(I_t)_{t=1, 2, \cdots}$ has no asymptotic influence on the existence and form of an equilibrium, and this is confirmed by comparing (35) and (36) with (30) and (31). 
Note however that the influence of $(I_t)$ can  be substantial on the development over time of effectives within the sub-populations. Furthermore, for any question of control needed to follow the society's obligation principle, the society
 has to look at each time $t$ at (36) and (37) so that all parameters governing $(I_t)$ as well as the integral $\Phi_{ni}(\tau_t)$ keep their relevance.
 
 \bigskip
  Sections 6, 7 and 8 constitute the main results of the present paper. The objective of the following Section is twofold. We want, on the one hand, advertise our approach in order to attract attention, and we do this by suggesting modifications of our model into different directions and showing that our approach stays compatible for several of them. On the other hand we would also like to briefly discuss how Society could envisage using claim distribution functions to effectively change the policy.
 
 \section{The flexibility of the approach}
\subsection{Modifications of our models}
\subsubsection{Modifying the attribution of ressources}
\noindent Our models inherit from the definition of a RDBP the
rule that individual claims are either served completely, or else, not at all, and that individuals with refused claims do not reproduce. One may object that this {\it zero-one} rule of satisfying claims is not fully realistic because individuals may accept compromises. It is relatively easy to modify our models into this direction without changing the essence our approach to find the new equilibrium criteria. 

Suppose that, as before, individuals of the sub-population $C_\delta, \delta \in \{h, i, ni\}$ submit  claims according to distribution $F_\delta$ and follow the law of reproduction $\{p_{\ell}^\delta\}_{\ell = 0, 1, 2, \cdots}, $ but Society follows the following scheme:
It decides to fix threshold values 
$$0=:s_0<s_1 <s_2<\cdots<s_k<s_{k+1}:=\infty,$$ say, and to partition offspring means in $k$ different classes:

\medskip \noindent
{\it Modification }:  An individual $j$ in class $C_\delta$ submitting the claim $X_j$ governed by the claim distribution $F_\delta$ will see its claim listed as
the maximum threshold below its claim, that is  $\tilde X_j= \sum_{l=1}^k s_l\one\{X_j\in [s_l,s_{l+1}[\}.$ Note that the applied policy does not depend on $\delta.$ If the claim $\tilde X_j$ is accepted we now suppose that the individual  will reproduce within $C_\delta$ with the conditional mean $$ \E\left(D^\delta(s_u)\Big| X_j \in [s_u, s_{u+1}[\,\right)=:m_\delta^u,~ \delta \in \{h, i, ni\}, u \in \{0, 1, \cdots k\} $$ being constant for each claim sub-interval. The idea is that the spacings can be chosen so that the desired effect is obtained, and that the probability of an individual being willing to reproduce can be built in these reproduction means.

For example, by choosing increasing gaps $s_{j+1}-s_j,$ the typical loss compared with the original claim will on the average increase with increasing claims, and a decreasing offspring mean over the claim intervals would reduce the willingness of an individual to stay and reproduce. But this can be expressed differently by modelling it in terms of an increasing probability $p_0^\delta$ of having no children, or again by a decresing mean reproduction. Choosing other spacings between the $s_j$'s
would have different effects. 

Now let us look at the effect on our results.

\medskip
It suffices to look at Model I. The random BRS-equation (13) as well as equation (15) only change in the sense that 

\smallskip  (i) the integrands of the integrals will change into $s_jdF_h(x)$, respectively $s_jdF_i(x)$ for $x\in [s_j,s_{j+1}[$ but $F_h$ and $F_i$ stay unchanged.

\smallskip (ii) the integrals now each split into a sum over the partition with different reproduction means which are constant over the sub-intervals of claims. 

\medskip \noindent Consequently, despite the truncation of the claims, the corresponding l.h.s. of (15) stays increasing and absolute continuous with respect to $\tau.$ The approach to find the equilibrium equation is thus the same. 

\medskip
\subsubsection{Modifying integration.}
There is also freedom for the definition of integration. If we suppose that that in generation $t$ the rate of integration is $\vfi_t$ with $0<\vfi_t<1$ then (27) and the random BRS-equation(28) = (29) keep their form with $\vfi$ being replaced by $\vfi_t.$ Thus, it suffices
to only assume that $(\vfi_t)_{t=1, 2, ...}$ converges to a limit $\vfi$ to maintain Theorem 7.1 as is. Of course,
the resulting limiting equilibrium $\alpha_\vfi$ may be quite different from $\alpha_\vfi$ for constant $\vfi.$
The "earlier history" evoked in 6.5.1 may have itself quite an erratic behavior if the sequence $(\vfi_t)$ converges slowly.

\medskip
Subsection 9.1 presents just two examples showing that the models can be generalised  without jeopardising the essence of our approach, and it is not hard to find other reasonable modifications for which our approach stays tractable.

\subsection{Changing from the {\it wf}-policy to more general policies}
\maketitle
So far we have always argued in terms of the counting random variable $N(n,s),$ respectively $\tilde N(n,s),$ which, as we recall here, is the relevant counting function for the weakest-first society ({\it wf}-society) (see Bruss and Duerinckx (2015), subsection 3.2.)
 The main reason is that according to Proposition 4.3 of this paper no society can survive
unless the {\it wf}-society can survive. 

The Theorem of envelopment states that the envelope of any development of a RDBP is formed by the {\it wf}-society and the {\it strongest first}-society ({\it sf}-society.) The latter has also a counting function based on order statistics, namely counting the decreasing order statistics of claims one can sum up without exceeding the available resource space. The lower bound given by the {\it sf}-society is not as neat as the upper bound, but it plays an equally important role because we know, again from the proof of the theorem of envelopment, that there can be no better one. 

The second reason is that the counting function $N(n,s),$ respectively $\tilde N(n,s),$ are in many ways the easiest to deal with. The idea to confine policies to order statistics of claims  is, as announced before, less restrictive than it seems, and we shall now extend on this fact. 

\subsection{ Order statistics based policies and extensions}

Let us  confine  our interest for the moment to a single RDBP with claim distribution $F.$
Trivially, no finite sum depends on the ordering of its terms,  so that any finite sum equals the sum of the increasing or decreasing order statistics of its terms. The same must hold for stopped sums provided that the quasi-stopping time defines the same selection of sum terms. The latter can often be obtained by adding additional rules leading to a de-facto change of the claim distribution function $F$.

We explain this by shortly discussing an example.

Suppose  Society may want to  favour middle-size claims at the expense of the lower range and/or upper range of claims. Society has many possibilities to achieve this.  For instance it can announce that, from the next generation onwards, claims exceeding a certain value $u$ (and only in this case), say will be submitted to a lottery $L$ with two possible outcomes $L\in\{0,1\}$. If $L=0$ then claim $X$ will be reduced to $cu$ for some fixed  $0< c\le 1$ and listed on the list of claims as $cu.$ If $L=1$ then it will be listed as $u+(X-u)/g(X)$ where $g(\cdot)$ is some fixed strictly increasing function with $g(\cdot)\geq 1$ everywhere. Let $\tilde X$ be the random outcome by submitting a claim $X.$ Hence
 $$\tilde X:= \one_{\{X\le u\}}\,X+\one_{\{X> u\}} \left(\one_{\{L=0\}}cu
 + \one_{\{L=1\}} \left\{u+(X-u)/g(X)\right\}\right).$$ With the decision, whether or not, to claim more than $u$ left to each individual, and $P(L=0)=a=1-P(L=1)$ each individual may look at the conditional expectation
\begin{align}\E(\tilde X|X) =\one_{\{X\le u\}}\,X +\one_{\{X> u\}}\big( acu+ (1-a)\{u+(X-u)/g(X)\}\big).\end{align}
There is nothing special about this new rule, the only fact we want to indicate here is to show how far one can already go with such a change of policy.
The constants $u, c, a$ as well as the function $g$ are chosen by the society which has thus an infinity of possibilities to push tendencies for claim sizes, in the desired direction. $\E(\tilde X|X)$ in  (44) shows a way to do this for  individuals who are supposed to rational with respect to expectation. Clearly, there are many other ways to influence tendencies.

The actual outcome of such interventions is not always easy to predict  because one would have to know
how those individuals, who would like to stay in the process, would relocate the claims they originally intended to submit. However, we know that if $\tilde F$ is the resulting new distribution of claims then the {\it wf}-policy for  $\tilde F$ corresponds to a different policy for $F,$ and to exemplify this was here are only modest objective. 

\smallskip 
The same reasoning holds of course for the claim-distributions $F_h, F_i,$ and $F_{ni}.$ Since any policy is supposed to be applied equally to all individuals, independently of to which sub-population they belong, this implies new constraints.
If we want to change $F_h$ into $\tilde F_h,$ and this change is formalised by an operator $\cal T,$ say, applied to $F_h,$ we have the constraint
\begin{align} \tilde F_h:= {\cal T}\circ F_h ~\implies \tilde F_i:={\cal T}\circ F_i {\rm~and~}\tilde F_{ni}:={\cal T}\circ F_{ni}\end{align} for the same operator $\cal T.$ It remains to propose tools we may use to construct such an operator $\cal T$ satisfying (46).

\subsection{Tools for controlling claim distributions} Although Society must apply the same policy to all individuals, it can nevertheless use all information it has about the socio-economic  structure of different layers of the sub-populations, such as e.g. unskilled workers, skilled workers, high-school levels, ... or different layers defined by a partition according to sex and age-classes, etc. If $F$ is the overall claim distribution function then Society sees $F$ as a mixture of distribution functions $F_j, ~j \in \{1, 2, \cdots, \ell\}$ where $\ell$ denotes the number of layers, and $F_j$ the distribution function of claims stemming from layer number $j.$

When Society considers how individuals in a sub-population, or within a certain layer of a sub-population may react to announcements of control such as tax policies, policies to encourage consumptions,  and others, then it sees a different kind of mixing. Any individual may change its claim behaviour by passing, independently of other individuals, with probability $g_j,$ say, to a  claim distribution $G_j.$  With additional information, Society may also consider a model assuming that an individual belonging to layer $i$
would change with probability $g^{(i)}_j,$ to a claim distribution $G^{(i)}_j.$

We see that there is a large number of possibilities to control claim distributions by different ways of mixing and we will not enter into questions of statistical nature which may naturally arise in practice. However, to keep the theory rigorous we now show that all our necessary conditions for equilibrium  stay rigorous, even if our "reference" distribution functions $F_h, F_i, F_{ni}$ are permanently  re-constituted over the generations by different kinds of mixing distributions. 

\subsection{BR-inequality, Steele's extension, and mixing distribution functions.}

We are here naturally interested in two kinds of mixing: On the one hand classical (independent) mixing where each random variable follows,  with a given probability $g_j$ a distribution function $F_j$, independently of the other random variables. On the other hand we have what one may call, in our context, {\it sub-class mixing} where an individual submits a $F_j$-distributed claim if and only if it belongs to some sub-class  of the {\it joint} population sharing the same characteristics 

\bigskip
\noindent {\bf Lemma 9~}
Let $k,n \in \N$ with $1\le k \le n$, and let $G$ be the mixed distribution function defined by $G(t)=\sum_{j=1}^k\,a_jF_j(t),$ where all $F_j$'s are absolute continuous distribution functions, and where all $a_j$'s are positive with $a_1+a_2+\cdots +a_k=1.$ Further let $Y_1, Y_2, \cdots, Y_n$ be  positive random variables with absolute continuous distribution function  $G.$ Let $N_Y(n,s)$ be the maximum number of them we can sum up without exceeding $s.$ If a value $\tau$ is the solution of  equation (9), 
i.e. $s=n\,\int_0^\tau x dG(x)$ then$$\E(N_Y(n,s))\le n\sum_{j=1}^k\,a_jF_j(\tau)$$ for both independent mixing and class mixing.

\medskip
\noindent {\bf Proof}: We first look at independent mixing for $n$ distributions. The absolute continuity of all $F_j(t)$ implies that $G(t),$ being a convex combination of the $F_j(t),$ is also absolute continuous. Hence, according to the BR-inequality (10),
\begin{align}\E(N_Y(n,s))\le n G(\tau) =n\sum_{j=1}^k\,a_jF_j(\tau).\end{align}

Second, for class mixing, let $X_1, X_2,  \cdots, X_n$ be a random sample of variables belonging to one of $1\le k\le n$ classes $C_j, j=1, \cdots, k$ with respective absolute continuous distribution functions $F_j.$ Let now $a_j=\#C_j/n $  and $G(x) = \sum_{j=1}^k a_jF_j(x).$ Let $X_{1,n}\le X_{2,n} \le \cdots \le X_{n,n}$ be the joint list of increasing order statistics and $\tilde N(n,s)$ be defined as in (7). Then we have \begin{align}~s\,=n\int_0^{\tau}x\,dG(x) ~\implies~ \E(\tilde N(n,s))\le \sum_{j=1}^k \#C_j \,F_j(\tau), \end{align}
since the values $s$ and $\tau$ on the l.h.s. satisfy (13) and $$s=n \int_0^\tau  x\, dG(x)= n \int_0^\tau  x\, d\Big(\sum_{j=1}^k a_j F_j(x)\Big)~~~~~$$$$~~~=\sum_{j=1}^k n a_j \int_0^\tau x\,dF_j(x) =\sum_{j=1}^k \#C_j\int_0^\tau x\,dF_j(x),$$ and since, with $a_j\ge 0$ and $\sum_{j=1}^k a_j=1$,   
 $G$ is again a convex combination of absolute continuous distribution functions. \qed
 
 \medskip
 \noindent {\bf Remark.} If we put $k=n$ and $\#C_j=1$ for all $j \in \{1,2,\cdots, n\}$ then this is Steele's extension. Since the BR-inequality is  a special case of Steele's extension, Lemma 9 is no more than a corollary of the latter. However, its interpretation has some benefit. Changing policies by re-defining  them as {\it wf}-policies for new mixed distributions will, whatever the mixing, and whatever the interdependence of mixing, even applied over successive generations, this will never  affect our bounds as long as we can trust on $F.$ As pointed out by Steele, the bounds only depend on the marginal distributions, and this means here very much.

\subsection{Controlling for convergence to an equilibrium}

Certain possibilities of control are easy to understand. If we look at the parameters $m_h, m_i,  r_h, r_i, ..$, and the important roles they play in our equations, they are prominent candidates for control to reach an equilibrium. However, in reality Society has often not so much influence on them. We have already mentioned deep-rooted cultural and traditional constraints to change birth rates, such as in particular $m_i.$ For other reasons it may be equally difficult
to control $r_h$ and $r_i.$ There may be  more freedom to change the rate $\vfi$ of integration, but both the home population and the immigrants have to cope with it, and there are certainly limits. Control must be sufficiently "smooth". 

Also, if a large $\vfi$ is in principle feasible, what would it mean for the home population if they become a minority in their country?
The evident factors are thus also the evident factors of control, but, in practice, they may be difficult to control, and it is clearly advisable to try to control several of them jointly (with smaller changes in mind)
than concentrating on one or two only.

\medskip What remains in our models are the claim distributions. Their influence is less evident, but for them Society will have many possible controls ranging from different tax policies over an uncountable set of possible initiatives up to target-oriented strict legislation. Whatever they might be, we only would like to discuss how to bring claim distributions into a desired direction. We want to comment on what guidelines Society may want to follow to do control efficiently, and in a sufficiently subtle way. This very last discussion is based on no more than a preliminary idea. Actually the author sees here a whole project of new research. 

\medskip\noindent
To make the subject of discussion precise it suffices confine our interest on model Model I. 
Recall that the conditions
\begin{align} m_hF_h(\tau)=m_iF_i(\tau)\ge1\end{align} are necessary for the existence of an equilibrium, where $\tau$ and $\alpha$ satisfy the limiting BRS-equality (15) with the parameters $r_h, r_i, m_h, m_i$ and the claim distribution functions $F_h$ and $F_i.$ 
Suppose that, for fixed $r_h, r_i, m_h, m_i$ there exists $b>0$ such that
$$\forall t\ge b:~ m_hF_h(t)\ge 1 \hbox{~and~} m_iF_i(t) \ge 1,$$ whereas $m_hF_h(s) \not= m_iF_i(s),$ for all $s \ge b.$ Then (49) cannot hold. Consequently,
the problem is now to find an operator ${\cal T}$ acting on the set of absolute continuous distribution functions such that \begin{align}m_h {\cal T}\circ F_h(\sigma) = m_i{\cal T}\circ F_i(\sigma)\end{align} for some $\sigma>b$ satisfying the corresponding limiting BRS-equation (15).
As we have seen Society can take measures which will will favour certain claims and discourage other claims and thus change the pre-supposed distribution function $F_h$ and $F_i$ into {\it de-facto} distribution functions
$\tilde F_h={\cal T}\circ F_h$ and $\tilde F_i={\cal T}\circ F_i$, say, so that claims will be identically distributed according to $\tilde F_h$
in the home-population and according to $\tilde F_i$
in the immigrant population.

\medskip
The following suggestion is no more that a preliminary idea:

We recall from measure theory that, given two measurable spaces $(S_1, \Sigma_1)$ and $(S_2, \Sigma_2)$, say, a measurable function $\psi:S_1\to S_2,$ and a measure $\mu: S_1 \to \R^+$, the so-called push-forward measure $\psi\#\mu: S_2\to \R^+$ is defined by \begin{align}
\forall A\in S_1: \psi \# \mu(A)=\mu(\psi^{-1}(A)).\end{align}
Concerning any question of optimality we recall a result for optimal transformation (transportation) on the real positive half-line (see e.g. Rachev and Rüschendorf (1998)): If the transport cost-function $c(x,y)$ measuring the cost per unit from $x$ to $y$ is convex in $|x-y|$ then there exists a unique optimal transport map realising the minimum, $c_{\rm min}$ say, namely
\begin{align} F_\nu^{-1}\circ\, F_\mu: \R\to \R ~\rm{and} ~c_{\rm min}=\int_0^1c(F_\nu^{-1}(t), F_\mu^{-1}(t))\,dt.\end{align}
Convexity in $|x-y|$ is a reasonable assumption, since costs (in our problem the cost of changing probability mass from $x$ to $y$) typically increase in a super-linear way.  For instance, we may choose   $c(x,y)=|x-y|^\beta$ for some $\beta>1.$

Now, to enable the existence of an equilibrium, we need not transform one  distribution function, $ F_h$ say, into some other specific distribution function $\tilde F_h,$ but, as seen in (49), we want to see them intersecting if weighted with $m_h$ and $m_i.$ Let us call this the "target."  

The target need not be reached  in one step. It can, at least in principle,  be planned in an {\it iterative} way. In a first step, we may apply class mixing and independent mixing 
to one side, $m_hF_h$ say, and check what this operation ${\cal T}_1$ does this to the other side. If the target is not reached, look which of the two resulting functions $$m_h{\cal T} _1\circ F_h~, ~m_i{\cal T} _1\circ F_i$$ would be more sensitive to control and mixing. We apply then this mixing to obtain two new sides $m_h{\cal T} _2 \circ{\cal T} _1\circ F_h$ and $m_i{\cal T} _2 \circ{\cal T} _1\circ F_i$, and so on. Few iteration  steps, (k, say) should do, and so, if methods of optimal transport can be applied stepwise, the hope is that
$${\cal T}:={\cal T}_k \circ {\cal T}_{k-1} \circ \cdots \circ {\cal T}_1$$ would be a close-to-optimal  pursuit of the target. \qed

\section{Critical comments and conclusions} 
\subsection{Critical comments}

1. Whatever their political or ideologic attitudes, most people agree that the hypotheses H1 and H2 are natural for human populations. However, H1, (populations want to survive forever), is based on a terminal event. Hence, provided that extinction is avoided in the next $n,$ say generations, the realisation of the terminal event depends only on Society's behaviour at times $n+1, n+1, \cdots.$  One only requires $P({\rm survival})>0,$ and no more. If the society obligation principle is ignored for a few generations, nothing essential may change. Does this mean that H1 should be better replaced by a more restrictive hypothesis?

The author's perception is as follows: Yes, with hypothesis H1, as it stands,  a limited number of exceptions to the society's obligation rule need not make an important difference.   The probability of survival will change, but it is true that the model is {\it under-determined} since only $P({\rm survival})>0$ is required. However, we 
propose to leave the hypothesis H1 as is because any more restrictive form having a true impact is likely to become too restrictive. This lies in the nature of the hypothesis. If parents insist that  their children have a future, would they care less for their grand-children, or grand-grand-children ...? If not, then, by recurrence,  a Society which does not respect its obligation principle in each generation, does not respect H1.
But then, what else could replace H1 as a generally accepted goal?

\medskip\noindent
2. Is the notion  of a {\it global model} as a sequence of unknown RDBPs not nebulous, or even redundant? Where is it used? 

\smallskip\noindent
No, it is neither nebulous nor redundant. RDBPs are well-defined, and so is the society's obligation principle. If Society examines in a generation the current RDBP and defines the next RDBP according to the principle (with all the freedom it has in doing so) then the current RDBP will have a successor. Thus we have a model, and recursively, a global model. Not knowing the specifications of the future RDBPs beforehand implies that we have no probability model. This comes with our intention to keep the model realistic but does not impair the character of being a model. We use the global model in Def. 6.2, Def. 7.1 and Def 8.2. which make no sense in local models.

\medskip \noindent
3. Where is the hypothesis H2 (standard of living) actually used?

\smallskip
\noindent
It is true that H2 never enters our mathematical analysis. However, it  intervenes in the freedom of choice of a policy in any generation. 
Bruss and Duerinckx (2015) (see Prop. 4.3, p. 346 and sect.7, pp 364-65) showed that the largest probability for survival and smallest expected standard of living (capacity to consume) per individual is obtained by the {\it wf}-policy. We cannot predict  how much Society would value, in any generation, a higher standard of living (by changing to another policy) on the expense of the survival probability. However,  under the classical hypothesis of scarce resources,  there is an upper bound for an affordable standard of living  under a given survival probability. 

\medskip\noindent
4. It is not realistic to suppose that the effectives of a human population can go to infinity. Therefore, what is the interest of criteria for the existence of equilibria conditioned on survival, i.e. conditioned on effectives tending to infinity?

\noindent The answer (linked with the one given in 3.) is that the asymptotical equilibrium is an idealisation showing us the way to understand necessary control in the view of H1 and H2 . The hypothesis H2 implies, as seen in Bruss and Duerinckx (2015), that Society will have the preference to keep the effective reproduction rate close to the critical value 1. Succeeding in doing so for a long time is what is really intended; unlimited growth is no convincing target for humanity.

\subsection{Conclusions}
Immigration leads to a large complex of different questions. Allowing for immigration is, on one side, an act of altruism, a grandeur of humanism.  On the other side, however, allowing for immigration may equally well be driven by lower motives, ranging in the worst case down to the intention to exploit a weaker sub-population.   The present article does not try to evaluate and compare advantages and disadvantages of immigration but only to study the question "When can it work out in the long run?"
The philosophy behind {\it working out} is that the inviting home-population may want to to keep essential parts  of its national identity, and the immigrants may want to import and live their own culture as far as these are not incompatible. This means in particular, none should wipe out the other one.  A long-term equilibrium between the respective effectives is seen as a necessary condition to make this possible,  and here it is Mathematics rather than Economics which must deliver answers. 

Although the mathematical arguments in this paper are rather elementary, tailoring our models in such a way that we can apply mathematics, and such that the models do not become unrealistic, this was more demanding. As pointed out, Steele's extension had much influence on our approach.

We have not discussed in detail the numerical computation of solutions, but, with the crucial functions $F_\delta (\tau)$ and $\Phi_\delta (\tau), \delta \in \{h,i,ni\}$ all being increasing in $\tau$, we had no difficulties in computing the solutions, as checked in several examples. Moreover, our results seem interesting. The mentioned sensitivity of the solutions as a function of the parameters is sometimes truly surprising, and it is good to have now explicit equilibrium criteria to see the exact reasons why. 

\smallskip
As far as the author is aware, the presented  approach to understand the mathematics of immigration and integration is new, and nothing comparable to the obtained explicit results has been known before. 
\subsection{Outlook} It would be convincing to  see our approach attract the attention of  specialists in optimal control/transport, and the interest of economists. We have provided the theory by taking care to justify all hypotheses. We have also tried our very best with respect to the transparency of the models and of the results they yield. However, for {\it reaching} the goals in the real life of  immigration, "smoothness" of actions must be expected to be equally important, and for this part our Subsection 9.6 may give no more than a modest beginning.\subsection*{Acknowledgement.} 

\smallskip\noindent
The author thanks M. Duerinckx for interesting discussions about pure-order policies, and L. R\"uschendorf for his expertise on  optimal transport problems.

\subsection*{ References}
~~~Afanasyev V. I., Geiger J.,  Kersting G. and V. A. Vatutin (2005) {\it Criticality for Branching Processes in Random Environment},
Ann.  Probability,
Vol. 33, No. 2,  645-673. 

\smallskip
Arlotto A., Mossel E. and Steele J.M. (2016)~ {\it Quickest online-selection of an increasing subsequence of specified size},  Random Structures \& Algorithms, Vol. 49, Issue 2, 235-252.

\smallskip Asmussen S. and Kurtz T. G.(1980)
{\it ~Necessary and Sufficient Conditions for Complete Convergence in the Law of Large Numbers},
       Ann. Probab.
    Vol. 8, Number 1, 176-182.

\smallskip Bansaye V., Caballero M.-E. and M\'el\'eard S. (2018){~\it Scaling limits of discrete population models with random environment and interactions},  4th Workshop on branch. processes, April 2018, Badajoz, Spain, Book of Abstracts, p. 26.

\smallskip
Barbour A., Hamza K.,  Kaspi H. and Klebaner F.C. (2015) ~{\it Escape from the boundary in Markov population processes}, Adv.  Appl. Prob., Vol 47, no. 4, 1190-1211. 

\smallskip
Barczy M., B\"osze Z. and Pap G. (2018) {\it Regularly varying non-stationary Galton-Watson processes with immigration}, arxiv.org/abs/1801.04002.

\smallskip
Bingham N. H. and R. A. Doney (1974) ~{\it Asymptotic properties of supercritical branching processes. I},  Adv.  Appl. Prob.,  Vol 6, 711-731.

\smallskip
Borjas, G. J. (2014) ~{\it Immigration Economics}, Harvard University Press, Cambridge Massachusetts, London.

\smallskip
Bruss F. T. (1978) ~{\it Branching Processes with Random Absorbing Processes}, J. Appl. Prob., Vol. 15, 54-64.

\smallskip
Bruss F. T. (1980) ~{\it A Counterpart of the Borel-Cantelli Lemma}, J. Appl. Prob., Vol. 17, 1094-1101.

\smallskip
Bruss F. T.  (1984a) ~{\it Resource Dependent Branching Processes}, 11th Conf. on Stoch. Proc. Applic. 1982, Abstract, Stoch. Proc. Applic., Vol. 16 (1), p. 36.

\smallskip
Bruss F. T.  (1984b) ~{\it A Note on Extinction Criteria for Bisexual Galton-Watson Processes}, J. Appl. Prob., Vol. 21, 915-919.

\smallskip
Bruss F. T. (2016) ~{\it The Theorem of Envelopment and Directives of Control in Resource
Dependent Branching Processes}, in Springer Lecture Notes in Statistics (I.M. del Puerto et al., Eds), Vol. 219, 119-136.

\smallskip
Bruss F. T. and Robertson J. B. (1991) ~{\it 'Wald's Lemma' for Sums of Order Statistics of i.i.d. Random Variables},
Adv.  Appl. Prob., Vol 23, 612-623.

\smallskip
Bruss F. T. and Delbaen F. (2001)~ {\it Optimal rules for the sequential selection of monotone subsequences of maximal length}, Stoch. Proc.  Applic., Vol 96, 313-342.

\smallskip
Bruss F. T. and Duerinckx M. (2015) ~{\it Resource dependent branching processes and the envelope of societies},
 Ann. of Appl. Probab., Vol. 25, Nr 1, 324-372.
 
\smallskip Daley D. J. (1968){ \it Extinction conditions for certain bisexual Galton-Watson branching processes}, Z. Wahrscheinlichkeitsth. und Verw. Geb.,
 Vol. 9, Issue 4,  315–322. 
 
  \smallskip 
 Gonz\'alez M., del Puerto I., and Yanev G., (2017) ~
 {\it Controlled Branching Processes}, J. Wiley \& Sons, Volume 2.

Haccou P.,  Jagers P. and  Vatutin V.A. (2005) ~{\it Branching Processes - Variation, Growth and Extinction of Populations}, Cambridge University Press.

\smallskip 
Hall P. and Heyde C.C. (1980)~ {\it Martingale Limit Theory and Its Applications}, Academic Press, New York, London, Toronto, Sydney, San Francisco.

\smallskip
Isp\'any M., Pap G. and van Zuijlen M.  (2005) ~{\it Fluctuation limit of branching processes with immigration and estimation of the means}, Adv. Appl. Probab. Vol 37, 523-538.

\smallskip  Jacquemain  A.  (2017){\it~Lorenz curves interpretations of the Bruss-Duerinckx Theorem for resource dependent branching processes}, arXiv:1708.01085v1.

\smallskip
Jagers, P. and Klebaner, F.C.  (2000) ~{\it Population-size-dependent and age-dependent branching processes}, Stoch. Proc. and Th. Applic., Vol 87,  235-254.

 \smallskip Keller G.,  Kersting  G. and R\"osler U. (1987){\it~ On the asymptotic behaviour of discrete time stochastic growth processes}, Ann. Probab. , Vol. 15, 305-343.

\smallskip
Kersting G. (2017), {~\it A unifying approach to branching processes
in varying environment},  ~arXiv:1703.01960v6 .

\smallskip 
Klebaner F.C. and  Zeitouni O. (1994) ~{\it The Exit Problem for a Class of Density-Dependent Branching Systems },
 Ann. of Appl. Probab.
Vol. 4, No. 4,  1188-1205. 

\smallskip
Molina M. (2010) ~{\it Two-sex branching process literature}, in Lecture Notes in Statistics,  196, Springer-Verlag 
(Eds. M. Gonzales, I.M. del Puerto, R. Martinez, M. Molina, M. Mota, A. Ramos), 279-291.

\smallskip Pap G. (2018), {\it Regularly varying Galton-Watson processes with immigration},  4th Workshop of Branch. proc. and their applic., April 2018, Badajoz, Spain, Book of Abstracts, p.\,71.

\smallskip
Perl I., Sen A. and Yadin A. (2015) ~{\it Extinction window of mean-field branching annihilating random walk,}
 Ann. of Appl. Probab., Vol. 25, Nr 6, 3139-3161.

\smallskip
Rachev S. and R\"uschendorf L. (1998) ~{\it Mass Transportation Problems}, 
Volume I, Theory, Springer-Verlag, New York, Berlin, Heidelberg.

\smallskip
Samuels S.M. and Steele J.M. (1981) ~{\it Optimal selection of a monotone subsequence from a random sample}, Ann. Probability,Vol. 9 (6), 937-947.

\smallskip
Steele J.M. (2016) ~{\it The Bruss-Robertson Inequality: Elaborations, Extensions, and Applications}, Math. Applicanda, Vol. 44,  No 1, 3-16.

\smallskip
Wajnberg A. (2014) ~{\it Le th\'eor\`eme de Bruss-Duerinckx ou l'enveloppement des soci\'et\'es humaines}, FNRS-News, Vol. 2014-06, pp (i)-(ii) and pp 20–22.

\smallskip
Wei C.Z. and Winnicki J.  (1989) ~{\it Some asymptotic results for the branching process with immigration}, Stoch. Proc. and Th. Applic., Vol. 31, 261-282.

\bigskip\bigskip

Author's address:
\smallskip

F.Thomas Bruss\hfill

Universit\'e Libre de Bruxelles,
Facult\'e des sciences, \hfill

D\'epartement de Mathématique, CP 210\hfill\

B-1050 Brussels, Belgium\hfill

tbruss@ulb.ac.be; phone:\,++32 2 650 5893.\hfill

\end{document}